
\documentstyle{amsppt}
\baselineskip18pt
\magnification=\magstep1
\pagewidth{30pc}
\pageheight{45pc}
\hyphenation{co-deter-min-ant co-deter-min-ants pa-ra-met-rised
pre-print pro-pa-gat-ing pro-pa-gate
fel-low-ship Cox-et-er dis-trib-ut-ive}
\def\leaderfill{\leaders\hbox to 1em{\hss.\hss}\hfill}

\def\Be{{\Cal B}}
\def\C{{\Cal C}}

\def\L{{\Cal L}}

\

\def\sgn{{\text {\rm \, sgn}}}
\def\idest{i.e.,\ }

\def\a{{\alpha}}
\def\be{{\beta}}
\def\g{{\gamma}}

\def\d{{\delta}}

\def\e{{\varepsilon}}
\def\z{{\zeta}}
\def\th{{\theta}}

\def\l{{\lambda}}

\def\cha{{\chi}}

\def\kase{{\tilde{e}}}
\def\kasf{{\tilde{f}}}

\def\b0{\text{\bf 0}}

\def\ra{{\ \longrightarrow \ }}

\def\lcon{{\leftarrow}}
\def\rcon{{\rightarrow}}

\def\hei{\text{\rm \, ht}}

\def\lan{{\langle}}
\def\ran{{\rangle}}

\def\lrt#1#2{\left\langle {#1}, {#2} \right\rangle}

\def\fg{{\frak g}}

\def\U{{U}}

\def\complex{{\Bbb C}}
\def\zed{{\Bbb Z}}
\def\kyu{{\Bbb Q}}

\def\Hom{\text{\rm Hom}}

\def\boxit#1{\vbox{\hrule\hbox{\vrule \kern3pt
\vbox{\kern3pt\hbox{#1}\kern3pt}\kern3pt\vrule}\hrule}}
\def\rabbit{\vbox{\hbox{\kern0pt
\vbox{\kern0pt{\hbox{---}}\kern3.5pt}}}}
\def\qchoose#1#2{\left[{{#1} \atop {#2}}\right]}

\def\tableau#1{
        \hbox {
                \hskip -10pt plus0pt minus0pt
                \raise\baselineskip\hbox{
                \offinterlineskip
                \hbox{#1}}
                \hskip0.25em
        }
}

\def\tabCol#1{
\hbox{\vtop{\hrule
\halign{\strut\vrule\hskip0.5em##\hskip0.5em\hfill\vrule\cr\lower0pt
\hbox\bgroup$#1$\egroup \cr}
\hrule
} } \hskip -10.5pt plus0pt minus0pt}

\def\CR{
        $\egroup\cr
        \noalign{\hrule}
        \lower0pt\hbox\bgroup$
}

\def\mapright#1{
    \mathop{\longrightarrow}\limits^{#1}}

\def\CD#1{
   $$
   \def\normalbaselines{\baselineskip20pt \lineskip3pt \lineskiplimit3pt }
    \matrix #1 \endmatrix
   $$
}

\def\blank#1#2{
\hbox to #1{\hfill \vbox to #2{\vfill}}
}


\def\strut{\vrule height10pt depth5pt width0pt}

\def\secta{1}
\def\sectb{2}
\def\sectc{3}
\def\sectd{4}
\def\secte{5}
\def\sectf{6}
\def\sectg{7}
\def\secth{8}
\def\secti{9}
\def\apage#1{}

\topmatter
\title Full heaps and representations of affine Kac--Moody 
algebras
\endtitle

\author R.M. Green \endauthor
\affil Department of Mathematics \\ University of Colorado \\
Campus Box 395 \\ Boulder, CO  80309-0395 \\ USA \\ {\it  E-mail:}
rmg\@euclid.colorado.edu \\
\newline
\endaffil

\abstract 
We give a combinatorial construction, not involving a presentation,
of almost all untwisted affine Kac--Moody algebras modulo their 
one-dimensional centres in terms of signed raising and lowering operators
on a certain distributive lattice $\Be$.
The lattice $\Be$ is constructed combinatorially as a set of ideals of 
a ``full heap'' over the Dynkin diagram, which leads to 
a kind of categorification of the positive roots for the Kac--Moody algebra.  
The lattice $\Be$ is also a crystal in the sense of Kashiwara, and its
span affords representations of the associated quantum affine algebra 
and affine Weyl group.
There are analogues of these results for two infinite families of twisted
affine Kac--Moody algebras, which we hope to treat more fully elsewhere.

By restriction, we obtain combinatorial 
constructions of the finite dimensional simple Lie algebras over 
$\complex$, except those of types $E_8$, $F_4$ and $G_2$.
The Chevalley basis corresponding to an arbitrary orientation of the Dynkin 
diagram is then represented explicitly by raising and lowering operators.
We also obtain combinatorial
constructions of the spin modules for Lie algebras of types $B$ and $D$,
which avoid Clifford algebras, and in which the action of Chevalley bases on
the canonical bases of the modules may be explicitly calculated.
\endabstract

 \endkeywords

\toc
\head{} Introduction \apage{2}\endhead
\head{1.} Heaps over Dynkin diagrams \apage{5}\endhead
\head{2.} Ideals of full heaps \apage{9}\endhead
\head{3.} Lie algebras and root systems \apage{16}\endhead
\head{4.} Parity of heaps in the simply laced case \apage{21}\endhead
\head{5.} Representability of roots in the simply laced finite type case 
 \apage{23}\endhead
\head{6.} The non simply laced case \apage{28}\endhead
\head{7.} Loop algebras and periodic heaps \apage{37}\endhead
\head{8.} Quantum affine algebras, crystals, and the Weyl group action
 \apage{43}\endhead
\head{9.} Applications and questions \apage{50}\endhead
\head{} Acknowledgements \apage{52}\endhead
\head{} Appendix: Examples of simply folded full heaps \apage{52}\endhead
\head{} References \apage{62}\endhead
\endtoc
\endtopmatter

\centerline{\bf To appear in the International Electronic Journal of Algebra}

\head Introduction \endhead

A heap is an isomorphism class of labelled posets, depending on an underlying 
graph $\Gamma$ and satisfying certain axioms.  Heaps have a wide variety of 
applications in algebraic combinatorics and statistical mechanics,
as explained in \cite{{\bf 21}}.  The algebraic and combinatorial theory of heaps
mostly concentrates on the case of finite heaps, but there is a well-developed
theory of infinite heaps used in the study of parallelism in computer
science, where they are known as ``dependence graphs'' \cite{{\bf 5}}.

In this paper, we introduce and study some remarkable infinite (but locally
finite) heaps, which we call ``full heaps'', and which have some interesting
applications to algebraic Lie theory.  Let $\Be$ denote the set of 
nonempty proper ideals of a full heap $E$ (regarded as a poset) over a 
graph $\Gamma$.  Using the 
poset structure, we will define a family of raising and lowering operators on
the space $V_E$ spanned by $\Be$.  If the underlying graph, $\Gamma$, of 
the full heap is a doubly laced Dynkin diagram associated to a symmetrizable
Kac--Moody algebra (meaning that all the entries of 
the corresponding generalized Cartan matrix lie in the set 
$\{2, 0, -1, -2\}$) then we will show how the space $V_E$ naturally carries 
the structure of (a) a module for the Kac--Moody algebra corresponding to 
$\Gamma$ and (b) a module for the Weyl group corresponding to $\Gamma$.
Moreover, the Chevalley generators (in the Kac--Moody case) and the Coxeter
generators (in the Weyl group case) act on $V_E$ via extremely simple raising
and lowering operators applied to basis elements.

If $\Gamma$ corresponds to an untwisted affine Kac--Moody algebra $\fg$, 
the representation of $\fg$ on $V_E$ over $\complex$ has a small kernel, 
namely the one-dimensional centre.  If we restrict attention to the
corresponding finite dimensional 
simple Lie algebra over $\complex$, the representation
will of course be faithful, but one can be much more precise: it is possible
to construct the Chevalley basis arising from a given orientation of the Dynkin
diagram (see \cite{{\bf 12}, (7.8.5), (7.9.3)}) explicitly in terms of raising 
and lowering operators.

Raising and lowering operators are familiar in other combinatorial
models of Lie theory.  The most important of these include the Kashiwara 
operators on crystals \cite{{\bf 13}}, used in the approach of the Kyoto 
school, and Littelmann's path operators \cite{{\bf 16}, {\bf 17}}.
Another example of raising and lowering operators occurs in the context
of the down-up algebras of Benkart and Roby \cite{{\bf 3}}.

For the case of simply laced finite dimensional
simple Lie algebras over $\complex$ (excluding
$E_8$), a combinatorial construction for the Lie algebras by raising and
lowering operators on ideals of (finite) heaps has been 
described by Wildberger \cite{{\bf 22}}.  Unfortunately, that paper
contains no proof of its main result \cite{{\bf 22}, Theorem 4.1}, which is 
analogous to our Theorem \sectf.7, and to the best of our knowledge,
no proof exists.  The constructions we describe 
here are modified versions of Wildberger's, and have some advantages over
them (see \S\secti).
Wildberger has also succeeded in dealing with the simple Lie algebra of type
$G_2$ using raising and lowering operators \cite{{\bf 23}}, but the construction
is ad hoc and significantly different from those of \cite{{\bf 22}} or of this 
paper.

The constructions described above require almost no knowledge of Lie theory,
apart from the definition of a Lie algebra and the notion of a Dynkin diagram
(or, equivalently, a generalized Cartan matrix).
In particular, the definition of a full heap is purely combinatorial.
However, the proofs that the constructions work do use Lie theoretic concepts.

Our representations of Kac--Moody algebras exist whenever we have a full heap
over the appropriate Dynkin diagram.  This includes all untwisted affine 
Kac--Moody algebras
except three (types $E_8^{(1)}$, $F_4^{(1)}$ and $G_2^{(1)}$ in Kac's
notation), and also includes two families of twisted affine Kac--Moody 
algebras 
(types $A_{2l-1}^{(2)}$ and $D_{l+1}^{(2)}$).  The more complicated root
systems in the twisted case make analysis more difficult,
so we will concentrate almost entirely on the untwisted
case in this paper for reasons of space.

Although our methods do not work for all types, they apply remarkably
uniformly 
in the cases where they do work.  The representations $V_E$ behave like
affine analogues of the minuscule representations of the corresponding
simple Lie algebras; in particular, the three cases mentioned above where
full heaps do not exist correspond to the three simple Lie algebras ($E_8$,
$F_4$ and $G_2$) that have no minuscule representations (see \cite{{\bf 19}, \S2.2} 
for more details).

The representation of a Kac--Moody algebra on $V_E$ has a $q$-analogue,
namely an action of the quantum affine algebra.  Regarded in this way,
$V_E$ is an integrable module and 
$\Be$ is a crystal basis for $V_E$ in the sense of Kashiwara, although
it does not give an extremal weight crystal.

Although we do not emphasise this in the sequel, the labelled heaps over
a fixed graph can be made into a category in which the isomorphism classes
of objects are precisely the heaps.
The framework of this paper can be regarded as a kind of categorification
of the positive roots of (most) affine Kac--Moody algebras, in which a
positive root $\a$ corresponds to a nonempty collection $\L_\a$ of 
labelled heaps.  An element of $\L_\a$ is called a root heap of character $\a$.
Isomorphic 
labelled heaps have the same characters, but since the isomorphism class is 
not determined by the character, we do not have a categorification in the
strict sense of \cite{{\bf 1}}, but rather in the weaker sense in which Khovanov 
homology \cite{{\bf 15}} is a categorification of the Jones polynomial.  This has the
consequence that root heaps have invariants that are not invariants of the
underlying positive root, and the most important for our purposes in the
simply laced case is that
of the parity of a root heap.  This depends on an arbitrarily chosen 
orientation of the Dynkin diagram and, when decategorified correctly
(Lemma \sectd.4), produces the asymmetry functions of
\cite{{\bf 12}, (7.8.4)}.  We also give a decomposition of a root heap into
convex sub-root heaps that corresponds to expressing a positive root as a
sum of two positive roots (Corollary \secte.5).  In the simply laced case,
this decomposition is unique, which corresponds to the fact that the structure
constants for the corresponding Chevalley basis lie in the set $\{-1, 0, 1\}$.
Our procedure for treating non simply laced cases is also a categorification
of a well-known procedure for producing non simply laced root systems, as
we discuss in \S\sectf.

The main results of the paper are as follows.  The representation of the
derived algebra $\fg'(A)$ of a symmetrizable Kac--Moody algebra is constructed
in Theorem \sectc.1.  The Chevalley bases for simple Lie algebras over
$\complex$ are constructed combinatorially in Theorem \sectf.7, and the
corresponding result for the whole untwisted affine Kac--Moody algebra modulo
its centre is given in Theorem \sectg.10.  A $q$-analogue of the latter
result is given in Theorem \secth.3, which also explains the connection
with crystal bases.

\head \S\secta. Heaps over Dynkin diagrams \endhead

Let $A$ be an $n$ by $n$ matrix with integer entries.  Following
\cite{{\bf 12}, \S1.1}, we call $A$ a {\it generalized Cartan matrix} 
if it satisfies the conditions (a) $a_{ii} = 2$
for all $1 \leq i \leq n$, (b) $a_{ij} \leq 0$ for $i \ne j$ and (c)
$a_{ij} = 0 \Leftrightarrow a_{ji} = 0$.  In this paper, we will only
consider generalized Cartan matrices with entries in the set 
$\{2, 0, -1, -2\}$; such matrices are sometimes called {\it doubly laced}.  
If, furthermore, $A$ has no entries equal to $-2$, 
we will call $A$ {\it simply laced}.

The Dynkin diagram $\Gamma = \Gamma(A)$ associated to a generalized Cartan
matrix is a directed graph, possibly with multiple edges, and vertices
indexed (for now) by the integers $1$ up to $n$.  If $i \ne j$ and
$|a_{ij}| \geq |a_{ji}|$, 
we connect the vertices corresponding to $i$ and $j$ by $|a_{ij}|$ lines;
this set of lines is equipped with an arrow 
pointing towards $i$ if $|a_{ij}|
> 1$.  For example, if $a_{ij} = a_{ji} = -2$, this will result in a double
edge between $i$ and $j$ equipped with an arrow pointing in each direction
(see Figure 7 in the Appendix).
There are further rules if $a_{ij} a_{ji} > 4$, but we do not need
these for our purposes.

The Dynkin diagram (together with the 
enumeration of its vertices) and the generalized Cartan matrix determine 
each other, so we may write $A = A(\Gamma)$.  If $\Gamma$ is 
connected, we call $A$ {\it indecomposable}.

Let $\Gamma$ be a Dynkin diagram with vertex set $P$ and no 
multiple edges.  Let $C$ be the relation on $P$ such that $x \ C\ y$ if and
only if $x$ and $y$ are distinct unadjacent vertices in $\Gamma$, and let
$\C$ be the complementary relation.

\definition{Definition \secta.1}
A {\it labelled heap} over $\Gamma$ is a triple $(E, \leq, \e)$ 
where $(E, \leq)$ is a locally finite partially ordered set (in other words,
a poset all of whose intervals are finite) with order relation denoted
by $\leq$ and where $\e$ is a map $\e : E \ra P$ satisfying the following two
axioms. 

\item{1.}{For every $\a, \be \in E$ such that $\e(\a) \ \C \ \e(\be)$, 
$\a$ and $\be$ are comparable in the order $\leq$.}

\item{2.}{The order relation $\leq$ is the transitive closure of the
relation $\leq_\C$ such that for all $\a, \be \in E$, $\a \ \leq_\C \ \be$ 
if and only if both $\a \leq \be$ and $\e(\a) \ \C \ \e(\be)$.}
\enddefinition

We call $\e(\a)$ the {\it label} of $\a$.  In the sequel, we will sometimes
appeal to the fact that the partial order is the reflexive, transitive 
closure of the covering relations, because of the local finiteness condition.

\definition{Definition \secta.2}
Let $(E, \leq, \e)$ and $(E', \leq', \e')$ be two labelled 
heaps over $\Gamma$.  We say that $E$ and $E'$ are isomorphic (as labelled
posets) if there is a poset isomorphism $\phi : E \ra E'$ such that 
$\e = \e' \circ \phi$.

A {\it heap} over $\Gamma$ is an isomorphism class of labelled
heaps.  We denote the heap corresponding to the labelled heap 
$(E, \leq, \e)$ by $[E, \leq, \e]$.
\enddefinition

We will sometimes abuse language and speak of the underlying set of a heap,
when what is meant is the underlying set of one of its representatives.

It can be shown \cite{{\bf 21}, \S2} that the finite heaps over a graph have a well
defined monoid structure, induced by an operation $\circ$ on labelled
heaps which we now define.

\definition{Definition \secta.3}
Let $E = (E, \leq_E, \e)$ and $F = (F, \leq_F, \e')$ be two finite labelled
heaps over $\Gamma$.
We define the finite labelled heap $G = (G, \leq_G, \e'') = E \circ F$ over
$\Gamma$ as follows.

\item{1.}{The underlying set $G$ is the disjoint union of $E$ and $F$.}
\item{2.}{The labelling map $\e''$ is the unique map $\e'' : G \ra P$ whose
restriction to $E$ (respectively, $F$) is $\e$ (respectively, $\e'$).}
\item{3.}{The order relation $\leq_G$ is the transitive closure of the
relation ${\Cal R}$ on $G$, where $\a \ {\Cal R} \ \be$ if
and only if one of the following three conditions holds:
\item{(i)}{$\a, \be \in E$ and $\a \leq_E \be$;}
\item{(ii)}{$\a, \be \in F$ and $\a \leq_F \be$;}
\item{(iii)}{$\a \in E, \ \be \in F$ and $\e(\a) \ \C \ \e'(\be)$.}
}
\enddefinition

\definition{Definition \secta.4}
Let $(E, \leq, \e)$ be a labelled heap over 
$\Gamma$, and let $F$ a subset of $E$.  
Let $\e'$ be the restriction of $\e$ to $F$.  Let ${\Cal R}$ be the relation
defined on $F$ by $\a \ {\Cal R} \ \be$ if and only if $\a \leq \be$ and
$\e(\a) \ \C \ \e(\be)$.  Let $\leq'$ be the transitive closure of ${\Cal R}$.
Then $(F, \leq', \e')$ is a labelled heap over $\Gamma$.  The heap 
$[F, \leq', \e']$ is called a {\it subheap} of $[E, \leq, \e]$.  

If $E = (E, \leq, \e)$ is a labelled heap over $\Gamma$, then we define
the {\it dual labelled heap}, $E^*$ of $E$, to be the labelled heap 
$(E, \geq, \e)$.  (The notion of ``dual heap'' is defined analogously.)

If $F$ is convex as a subset of $E$ (in other words, if 
$\a \leq \be \leq \g$ with $\a, \g \in F$, then $\be \in F$) then we call 
$F$ a {\it convex subheap} of $E$.  If, whenever $\a \leq \be$ and $\be 
\in F$ we have $\a \in F$, then we call $F$ an {\it ideal} of $E$; dually,
if, whenever $a \geq \be$ and $\be \in F$ we have $\a \in F$, then we call
$F$ a {\it filter} of $E$.
If $F$ is an ideal of $E$ with $\emptyset \subsetneq F \subsetneq E$ 
such that for each vertex $p$ of $\Gamma$ we have 
$\emptyset \subsetneq F \cap \e^{-1}(p) \subsetneq \e^{-1}(p)$, then we
call $F$ a {\it proper ideal} of $E$.
\enddefinition

\remark{Remark \secta.5}
If $E$ and $F$ are finite heaps over $\Gamma$, then it follows from the above
two definitions that $E$ and $F$ are both convex subheaps of $E \circ F$,
and that $E$ is an ideal of $E \circ F$.
\endremark

We will often implicitly use the fact that a subheap is determined by its
set of vertices and the heap it comes from.
 
\definition{Definition \secta.6}
Let $(E, \leq, \e)$ be a locally finite labelled heap over $\Gamma$.
We say that $(E, \leq, \e)$ and $[E, \leq, \e]$ are {\it fibred}
if
\item{(a)}{for each vertex $p$ in $\Gamma$, the subheap $\e^{-1}(p)$ is
unbounded above and unbounded below,}
\item{(b)}{for every pair $p, p'$ of adjacent vertices in $\Gamma$
and every element $\a \in E$ with $\e(\a) = p$, there exists $\be \in E$
with $\e(\be) = p'$ such that either $\a$ covers $\be$ or $\be$ covers $\a$
in $E$.}
\enddefinition

\remark{Remark \secta.7}
\item{(i)}{It is easily checked that these are sound definitions, because 
they are invariant under isomorphism of labelled heaps.}
\item{(ii)}{The name ``fibred'' alludes to the fact that these heaps can
also be constructed using fibre bundles.  For $x \in E$, define the set
$O_x^E \subseteq E \times \Gamma$ to consist of all pairs $(x, p)$, where
$x \in E$ and there exists $y \in E$ with $\e(y) = p$ such that
either $x$ covers $y$ or $y$ covers $x$.  For each vertex $a$ of $\Gamma$, 
define the set $O_a^\Gamma \subseteq \Gamma \times \Gamma$ to consist of 
all pairs $(a, b)$ such that $a$ and $b$ are adjacent in $\Gamma$.  Define
$\pi : E \times \Gamma \ra \Gamma \times \Gamma$ by $\pi((x, p)) = (\e(x), p)$.
Let $\zed$ be the set of integers equipped with the discrete topology,
equip $E \times \Gamma$ with the smallest topology such that the sets $O_x^E$
are open, and equip $\Gamma \times \Gamma$ with the smallest topology such
that the sets $O_a^\Gamma$ are open.  Then $E$ is fibred if and only if
\CD{
\zed & \mapright{} & E \times \Gamma & \mapright{\pi} & \Gamma \times \Gamma
}
is a fibre bundle.
}
\item{(iii)}{Condition (a) provides a way to name the elements of $E$, which
we shall need in the sequel.
Choose a vertex $p$ of $\Gamma$.  Since $E$ is locally finite, $\e^{-1}(p)$
is a chain of $E$ isomorphic as a partially ordered set to the integers, so
one can label each element of this chain as $E(p, z)$ for some $z \in \zed$.
Adopting the convention that $E(p, x) < E(p, y)$ if $x < y$, this labelling
is unique once a distinguished vertex $E(p, 0) \in \e^{-1}(p)$ has been
chosen for each $p$.}

\definition{Definition \secta.8}
Let $E$ be a fibred heap over a Dynkin diagram $\Gamma$ with generalized
Cartan matrix $A$.  If every open interval $(\a, \be)$ of $E$
such that $\e(\a) = \e(\be) = p$ and $(\a, \be) \cap \e^{-1}(p) = \emptyset$
satisfies $\sum_{\g \in (\a, \be)} a_{p, \e(\g)} = -2$, we call $E$ a
{\it full} heap.
\enddefinition

The above definition is reminiscent of Stembridge's definition of a 
minuscule heap in \cite{{\bf 20}, \S3}; however, we are following Kac's 
definition of generalized Cartan matrix \cite{{\bf 12}, \S4.7}, which is the 
transpose of Stembridge's.  (This distinction only applies to the
matrices, and not to the corresponding heaps.)  The definition says that 
either (a) $(\a, \be)$ contains precisely two elements with labels ($q_1$, 
$q_2$, say) adjacent to $p$ and such that there is no arrow from 
$q_1$ (or $q_2$) to $p$ in the Dynkin diagram, or that (b) $(\a, \be)$ 
contains precisely one element with label ($q$, say) adjacent to $p$ 
such that there is an arrow from $q$ to $p$ in the Dynkin diagram.

\head \S\sectb. Ideals of full heaps \endhead

In \S\sectb, we develop some properties of ideals of full heaps, and use them
to define raising and lowering operators.

\proclaim{Lemma \sectb.1}
Let $(E, \leq, \e)$ be a full labelled heap over 
$\Gamma$, let $F$ and $F'$ be proper ideals of $E$ and let $J$ be an
ideal of $E$.
\item{\rm (i)}{With the labelling convention of Remark \secta.7,
the ideal $J$ is proper if and only if for all vertices $p \in \Gamma$,
we have $$
J \cap \e^{-1}(p) = E_p(N) := \{ E(p, t) : t < N \} 
$$ for some integer $N$ depending on $J$ and $p$.}
\item{\rm (ii)}{The subheaps $F \cap F'$ and $F \cup F'$ are proper ideals
of $E$.}
\item{\rm (iii)}{If $F \subseteq J$ and $J \backslash F$ is finite, then $J$
is a proper ideal of $E$.}
\item{\rm (iv)}{If $J \subseteq F$ and $F \backslash J$ is finite, then $J$
is a proper ideal of $E$.}
\item{\rm (v)}{If $\Gamma$ is finite and connected, then $J$ is a 
proper ideal if and only if $\emptyset \ne J \ne E$.}
\item{\rm (vi)}{If $\Gamma$ is finite and $F \subseteq F'$,
then $F' \backslash F$ is finite.}
\item{\rm (vii)}{If $\Gamma$ is finite and connected and $L$ is a finite
convex subheap of $E$, then there is a proper ideal $L'$ of $E$ such that
$L \subset L'$ and $L' \backslash L$ is a proper ideal of $E$.}
\endproclaim

\demo{Proof}
Part (i) follows from Remark \secta.7 (iii) and the definition of proper 
ideal.

Part (ii) follows from (i) and the fact that, for a fixed vertex $p$, the 
chains $E_p(N)$ are closed under finite intersections and finite unions.

For part (iii), choose $N$ so that $F \cap \e^{-1}(p) = E_p(N)$.  Since
$J \backslash F$ is finite, $$
(J \cap \e^{-1}(p)) \backslash  (F \cap \e^{-1}(p))
$$ must be finite.  However, since $J$ is an ideal,  
$J \cap \e^{-1}(p)$ must be downward closed (meaning that if $E(p, y)
\in J \cap \e^{-1}(p)$ and $x < y$, then $E(p, x) \in J \cap \e^{-1}(p)$).
This shows that $J \cap \e^{-1}(p) = E_p(N')$ for some $N' \geq N$.

Part (iv) follows by a similar argument to that used to prove part (iii),
mutatis mutandis.

If $J$ is a proper ideal, then it follows easily from (i) that $\emptyset
\ne J \ne E$, as required for (v), 
so suppose that $\emptyset \ne J \ne E$ for an ideal $J$.
Let $\a \in E \backslash J$ and let $p = \e(\a)$.

Suppose first that $J \cap \e^{-1}(p) \ne \emptyset$.  Since 
$J \cap \e^{-1}(p)$ is an ideal of $\e^{-1}(p)$ and $\a \not\in J$, we have
$J \cap \e^{-1}(p) = E_p(N)$.  Let $q$ be adjacent to $p$ in $\Gamma$.
The definition of full heap ensures that there exists an element $\be > \a$
with $\e(\be) = q$, and since $\a \not \in J$, we must have $\be \not\in J$.
On the other hand, there also exists an element $\be' < E(p, N-1)$ with
$\e(\be') = q$, and the fact that $J$ is an ideal means that $\be' \in J$.
Combining these observations, we see that $J \cap \e^{-1}(q) = E_q(N')$ for
some integer $N'$.  Since $\Gamma$ is connected, a similar condition holds
at each vertex, and $J$ is proper.

The other possibility is that $J \cap \e^{-1}(p) = \emptyset$.  By reversing
the argument of the above paragraph, we find that $J \cap \e^{-1}(q) = 
\emptyset$ for all vertices $q$, in other words, that $J = \emptyset$,
contrary to the hypothesis of (v).

Under the assumptions of (vi), the sets $(F' \backslash F) \cap \e^{-1}(p)$ are
all finite by (i), which means that $F' \backslash F$ is also finite because
$\Gamma$ is.

For part (vii), define $L' = \{\a \in E : \a \leq \be \text{ for some }
\be \in L \}$.  It is easily checked that $L'$ is a nonempty ideal of $E$ and
that $L' \backslash L$ is an ideal of $E$.
Furthermore, $L'$ is bounded above (because $L$ is), so $L' \ne E$, 
$L$ is a proper ideal by (v) and $L' \backslash L$ is a proper ideal
by (iii).  
\qed\enddemo

Part (vi) of Lemma \sectb.1 will often be used without comment in the
sequel.  Part (ii) of the lemma has the following immediate corollary.

\proclaim{Corollary \sectb.2}
The set of all proper ideals of $E$ has the structure of a distributive
lattice, where $I \wedge J := I \cap J$ and $I \vee J := 
I \cup J$.
\qed\endproclaim

\definition{Definition \sectb.3}
Let $R^+$ be the set of all functions $P \ra \zed^{\geq 0}$.
If $F$ is a finite labelled heap over $\Gamma$, then 
we define the {\it character}, $\cha(F)$ of $F$ to
be the element of $R^+$ such that $\cha(F)(p)$
is the number of elements of $F$ with $\e$-value $p$.  
If $\a \in R^+$, we write $\L_\a(E)$ to be the set of all convex
subheaps $F$ of $E$ with $\cha(F) = \a$.  If $F$ consists
of a single element $\a$ with $\e(\a) = p$, we will write $\cha(F) = p$ for
short, so that $\L_p(E)$ is identified with the elements of $E$ labelled
by $p$.

Since the function $\cha$ is an invariant of labelled heaps, we can extend
the definition to apply to finite heaps of $\Gamma$.
\enddefinition

\example{Example \sectb.4}
Let $\Gamma$ be the Dynkin diagram of type $E_7^{(1)}$, shown in Figure
17 in the Appendix, let $E$ be the heap shown in Figure 18, and let
$F$ be the finite convex subheap shown in the dashed box.  Writing $\a_i$ for
the function sending $i \in P$ to $1$ and $j \in P$ to $0$ if $j \ne i$,
we find that $$\cha(F) = \a_0 + 2 \a_1 + 3 \a_3 + 4 \a_3 + 3 \a_4 + 2 \a_5
+ \a_6 + 2 \a_7.$$
\endexample

\proclaim{Lemma \sectb.5}
Let $[E, \leq, \e]$ be a full heap over the Dynkin diagram $\Gamma$, with
generalized Cartan matrix $A$, let $I$ be an ideal of $E$ 
and let $\a \in I$ be a maximal element.  Define $p = \e(\a)$, and suppose 
that $q \in \Gamma$ is adjacent to $p$.  Then precisely one of the following 
occurs:
\item{\rm (i)}{$a_{qp} = -1$, $I \backslash \{\a\}$ is an ideal of $E$ and 
there exists a
maximal element $\be \in I \backslash \{\a\}$ with $\e(\be) = q$,}
\item{\rm (ii)}{$a_{qp} = -1$, there exists a minimal element 
$\be \in E \backslash I$ such 
that $\e(\be) = q$ and $I \cup \{\be\}$ is an ideal of $E$, or}
\item{\rm (iii)}{$a_{qp} = -2$, there exists a maximal element $\be \in I 
\backslash \{\a\}$ and a minimal element $\be' \in E \backslash I$ such that
$\e(\be) = \e(\be') = q$, and both $I \backslash \{\a\}$ and $I \cup \{\be'\}$
are ideals of $E$.}
\endproclaim

\demo{Proof}
By part (b) of the definition of a fibred heap, there exists $\be \in E$ with
$\e(\be) = q$ such that either $\a$ covers $\be$ or $\be$ covers $\a$.
Until further notice, let us assume that $a_{qp} = -1$.

Suppose first that $\be < \a$.  Since $\a$ is maximal in $I$, it follows
that $I \backslash \{\a\}$ is an ideal of $E$.  If $\be$ is maximal in
$I \backslash \{\a\}$, then we are in the situation of (i) above, so suppose
this is not the case.  By part (a) of the definition of a fibred heap, 
there exists $\g \in E$
with $\e(\g) = q$ and $\g > \be$.  Since $\{\a, \be, \g\}$ is a chain in $E$
and $\be < \a$ is a covering relation, we have $\g > \a$ and $\g \not\in I$.  
Because $\e^{-1}(q) \cup \{\a\}$ is a chain in $E$, we may assume that the 
interval
$(\be, \g)$ of $E$ contains no elements of $\e^{-1}(q)$.  By the definition
of full, $(\be, \g)$ contains two elements, $\g_1$ and $\g_2$, with labels 
adjacent to $q$, and one of these elements, $\g_1$ say, is $\a$.  The 
hypothesis that $\be$ is not maximal in $I \backslash \{\a\}$ implies that
$\g_2$ is also in $I$.  We claim now that $I \cup \{\g\}$ is an ideal of $E$;
to show this, it is enough to show that if $\g > \g'$ is a covering relation,
then $\g' \in I$.  The latter holds because any such $\g'$ is comparable to 
$\be$, which has the same label as $\g'$, and $I \cup \{\g\}$ contains all 
elements less than $\be$ together with the closed interval $[\be, \g]$.
This satisfies the conditions of (ii).

Suppose now that $\be > \a$.  It is clear that $\be$ is minimal in $E
\backslash I$, so that if $I \cup \{\be\}$ is an ideal of $E$, the conditions
of (ii) will hold.  Suppose that this is not the case.  By part (a) of the
definition of a fibred heap, 
there exists $\g \in E$ with $\e(\g) = q$ and $\g < \be$.  As in the
previous paragraph, this means that $\g < \a$, from which we see that
$\g \in I$.   By the definition of full, $(\g, \be)$ contains two elements, 
$\g_1$ and $\g_2$, with labels adjacent to $q$, and one of these elements, 
$\g_1$ say, is $\a$.  If $\g_2$ also lies in $I$, then the conditions of
(ii) will be satisfied as in the previous paragraph.  If $\g_2$ does not
lie in $I$, then $I \cup \{\be\}$ is not an ideal, but $I \backslash \{\a\}$
is an ideal with maximal element $\g$, and $\e(\g) = q$.  This satisfies
the conditions of (i).

From now on, assume that $a_{qp} = -2$.  Suppose there exists $\be \in E$
with $\e(\be) = q$ such that $\a$ covers $\be$.  By part (a) of the definition
of a fibred heap, there exists $\g \in E$
with $\e(\g) = q$ and $\g > \be$.  Following the same reasoning as earlier,
we have $\g > \a$ and $\g \not\in I$.
The other possibility is that there exists $\be \in E$
with $\e(\be) = q$ such that $\be$ covers $\a$.  In this case, part (a) of the
definition of a fibred heap shows that
there exists $\g \in E$ with $\e(\g) = q$ and $\g < \be$.
As before, this means that $\g < \a$, from which we see that
$\g \in I$.  In either case, there exists a chain $\be_1 < \a < \be_2$ in
$E$ with $\e(\be_1) = \e(\be_2) = q$, such that
the open interval $(\be_1, \be_2)$ contains no elements labelled $q$.
By the definition of full, $\a$ is the only element in $(\be_1, \be_2)$ with
a label adjacent to $q$, and thus the only element in $(\be_1, \be_2)$,
meaning that $\{\be_1, \a, \be_2\}$ is a convex chain.  The assertions of
case (iii) now follow by adapting the argument for the case $a_{qp} = -1$.
\qed\enddemo

The following definition generalizes ideas in \cite{Wi, \S4}.

\definition{Definition \sectb.6}
Let $E$ be a full heap over a graph $\Gamma$, let $k$ be a field (always
of characteristic not equal to $2$), and define 
$V_E$ to be the $k$-span of the distributive lattice $$
\Be = \{v_I : I \text{ is a proper ideal of } E \}
.$$  For any such ideal and any finite convex subheap $L \leq E$, we
write $L \succ I$ to mean that both $I \cup L$ is an ideal and $I \cap L = 
\emptyset$, and we write $L \prec I$ to mean
that both $L \leq I$ and $I \backslash L$ is an ideal.
We define linear operators 
$X_L$, $Y_L$ and
$H_L$ on $V_E$ as follows: $$\eqalign{
X_L(v_I) &= \cases
v_{I \cup L} & \text{ if } L \succ I,\cr
0 & \text{ otherwise,}\cr
\endcases\cr
Y_L(v_I) &= \cases
v_{I \backslash L} & \text{ if } L \prec I,\cr
0 & \text{ otherwise,}\cr
\endcases\cr
H_L(v_I) &= \cases
v_I & \text{ if } L \prec I \text{ and } L \not\succ I,\cr
-v_I & \text{ if } L \succ I \text{ and } L \not\prec I,\cr
0 & \text{ otherwise.}\cr
\endcases\cr
}$$  (Note that these operators are defined by parts (iii) and (iv) of 
Lemma \sectb.1; they are nonzero by part (vii) of Lemma \sectb.1.)
If $p$ is a vertex of $\Gamma$, we write $X_p$ for the linear operator
on $V_E$ given by $\sum_{L \in \L_p(E)} X_L$ (with notation as in Definition
\sectb.3), and we define $Y_p$ and $H_p$
similarly.  Note that although these sums are infinite, it follows from
the definitions of fibred and full heaps that at most one of the
terms in each case may act in a nonzero way on any given $v_I$.  In this
situation, we also write $p \succ I$ to mean that $L \succ I$ for some
(necessarily unique) $L \in \L_p(E)$, and analogously we write
$p \prec I$ with the obvious meaning.  Note that it is not possible for
both $p \prec I$ and $p \succ I$, because $I$ cannot contain a convex
chain $\a < \be$ with $\e(\a) = \e(\be) = p$.
\enddefinition

\proclaim{Lemma \sectb.7}
Maintain the above notation and suppose that $p$ and $q$ are vertices of
$\Gamma$ (allowing the possibility $p = q$).  We have the following relations 
in the associative $k$-algebra
generated by the operators $X_L$, $Y_L$ and $H_L$, where $\d$ is the
Kronecker delta: $$\eqalignno{
H_p H_q &= H_q H_p, & (1)\cr
H_p X_q - X_q H_p &= a_{pq} X_q, & (2)\cr
H_p Y_q - Y_q H_p &= -a_{pq} X_q, & (3)\cr
X_p Y_q - Y_q X_p &= \d_{pq} H_q, & (4)\cr
X_p X_q &= X_q X_p,  \text{ if } a_{pq} = 0, & (5)\cr
Y_p Y_q &= Y_q Y_p,  \text{ if } a_{pq} = 0, & (6)\cr
X_p X_p &= Y_p Y_p = 0, & (7)\cr
X_p X_q X_p &= Y_p Y_q Y_p = 0  \text{ if } a_{pq} = -1. & (8)\cr
}$$
\endproclaim

\demo{Proof}
Relation (1) holds because of the way the operators $H_p$ act as scalars on 
each basis vector $v_I$.

Consider the algebra element $H_p X_p - X_p H_p$.  This element, and $X_p$, 
will each act as zero on $v_I$ unless $p \succ I$.  If, on the other hand,
$p \succ I$, let $L \succ I$ be such that $L = \{\a\}$ with
$\e(\a) = p$.  We then have $
H_p X_p v_I = v_{I \backslash L} = - X_p H_p v_I
,$ and relation (2) follows.  

Suppose that $a_{pq} = 0$, so that $p$ and $q$ are not adjacent, and 
consider $H_p X_q - X_q H_p$.
Unless $q \succ I$, this element will act as zero on $v_I$, so we may reduce
consideration to this case.  
Let $L \succ I$ be such that $L = \{\a\}$ with
$\e(\a) = q$.  
A simple case by case check shows that $p \prec I \cup L$
(respectively, $p \succ I \cup L$) if and only $p \prec I$ (respectively,
$p \succ I$), and relation (2) follows.

Now suppose that $a_{pq} = -1$, and consider $H_p X_q - X_q H_p$.
As before, relation (2) follows trivially unless $q \succ I$, so we
reduce to this case.  Let $L \succ I$ be such that $L = \{\a\}$ with
$\e(\a) = q$, so that $X_q v_I = v_{I \cup L}$.  By Lemma \sectb.5, we have 
either $p \succ I \cup L$ or $p \prec I$, but not both.  (Note that if
$p \succ I \cup L$, then $p \not\succ I$, and if $p \prec I$, then $p 
\not\prec I \cup L$.)  If 
$p \succ I \cup L$, we have $H_p X_q v_I = - v_{I \cup L}$ and 
$X_q H_p v_I = 0$.  On the other hand if $p \prec I$, we have
$H_p X_q v_I = 0$ and $X_q H_p v_I = v_{I \cup L}$.  Relation (2) now follows.

Finally, suppose that $a_{pq} = -2$, and consider $H_p X_q - X_q H_p$.
Again, let $L \succ I$ be such that $L = \{\a\}$ with
$\e(\a) = q$, so that $X_q v_I = v_{I \cup L}$.
By Lemma \sectb.5, we have both $p \succ I \cup L$ and $p \prec I$,
meaning that $H_p X_q v_I = - v_{I \cup L}$ and $X_q H_p v_I = -v_{I \cup L}$.
This completes the proof of relation (2).

The verification of relation (3) follows a similar line of argument to that
used to prove relation (2), mutatis mutandis.

It follows from the definition of $H_p$ that $X_p Y_p - Y_p X_p = H_p$, thus
establishing the case $p = q$ of relation (4).  If $p$ and $q$ are adjacent, 
then a case by case check shows that the operators $X_p Y_q$ and $Y_q X_p$ are
individually zero.  On the other hand, if $p$ and $q$ are not adjacent,
an argument like that used on the $a_{pq} = 0$ case of relation (2) shows that
$X_p Y_q$ and $Y_q X_p$ commute, thus finishing the proof of relation
(4).  The same reasoning also establishes the commutation relations (5) and 
(6).

Relation (7) holds because no ideal $I$ can contain a convex chain $\a < \be$
with $\e(\a) = \e(\be) = p$, and relation (8) holds because no ideal $I$
can contain a convex chain $\a < \be < \g$ with $\e(\a) = \e(\g) = p$ and
$\e(\be) = q$ if $a_{pq} = -1$.  (Both of these conditions come from the 
definition of full.)
\qed\enddemo

\head \S\sectc.  Lie algebras and root systems \endhead

A {\it Lie algebra} over a field $k$ is a $k$-vector space $\fg$ endowed with a
bilinear (usually nonassociative) multiplication $[\, , \, ] : \fg \times \fg
\ra k$.  The image of the pair $(x, y)$ under this map is denoted by $[x, y]$,
and the following axioms hold for all elements $x, y, z \in \fg$: $$\eqalign{
[x, x] &= 0;\cr
[x, [y, z]] + [y, [z, x]] + [z, [x, y]] &= 0.\cr
}$$  The first condition above is known as {\it antisymmetry} and the second
is known as the {\it Jacobi identity}.  Any associative algebra $A$ over $k$, 
such as the algebra of Lemma \sectb.7,
may be made into a Lie algebra using the bracket $[a, b] := ab - ba$.

The significance of Lemma \sectb.7 is that it gives $V_E$ the structure
of a module for a certain Lie algebra, namely (in the case $k = \complex$)
the derived algebra $\fg'(A)$
of a symmetrizable Kac--Moody algebra (see \cite{{\bf 12}, \S0.3}).
The main purpose of this paper is to understand this module.

\proclaim{Theorem \sectc.1}
Let $E$ be a full heap over a Dynkin diagram $\Gamma$ with vertices $P$ and
generalized Cartan matrix $A$, let $k$ be 
a field of characteristic different from $2$, and let $V_E$ be the 
$k$-vector space of Definition \sectb.6.
Let $\fg$ be the Lie algebra with generators $\{e_i, f_i, h_i : i \in P\}$
and the usual defining relations (see \cite{{\bf 4}, \S9.4}).
Then $V_E$ becomes a left $\fg$-module, where $e_i$ (respectively,
$f_i$, $h_i$) acts as $X_i$ (respectively, $Y_i$, $H_i$).
\endproclaim

\demo{Proof}
This is a consequence of Lemma \sectb.7, recalling that we have $$
[a, b] . v := a(b(v)) - b(a(v))
.$$  Note that the relation $[e_i, [e_i, e_j]] = 0$ expands to $$
e_i \circ e_i \circ e_j - e_i \circ e_j \circ e_i + e_j \circ e_i \circ e_i = 0
,$$ which holds by relations (7) and (8) of Lemma \sectb.7, and the relation
$$[e_i, [e_i, [e_i, e_j]]] = 0$$ expands to $$
e_i \circ e_i \circ e_i \circ e_j 
- 2e_i \circ e_i \circ e_j \circ e_i
+ 2e_i \circ e_j \circ e_i \circ e_i 
- e_i \circ e_j \circ e_i \circ e_i 
= 0
,$$ which holds by relation (7) of Lemma \sectb.7.  Similar
comments hold for the relations involving $f_i$ and $f_j$.
\qed\enddemo

To understand Lie algebras such as those in Theorem \sectc.1, one needs the 
concept of a root system, and we will show that 
the combinatorics of full heaps is intimately connected to that of root
systems for affine Kac--Moody algebras.  

We introduce the following notation in order to state later results easily.

\definition{Definition \sectc.2}
Let $L$ be a finite convex heap of a full heap $E$ over a graph $\Gamma$, and
let $p$ be a vertex of $\Gamma$.  We write $p \rcon L$ (respectively, 
$L \lcon p$) to mean that $L$ has a minimal (respectively, maximal) vertex 
with label $p$.  We write $p \lcon L$ (respectively, $L \rcon p$) to mean
that there is a (necessarily unique) vertex $\a$ of $E \backslash L$ 
labelled $p$ such that $L \cup \{\a\}$ is convex and $p \rcon L \cup \{\a\}$ 
(respectively, $p \lcon L \cup \{\a\}$).  Define the integers $b^\pm(L, p)$
by the conditions $$
b^+(L, p) = \cases
1 & \text{ if } L \lcon p,\cr
-1 & \text{ if } L \rcon p,\cr
0 & \text{ otherwise},
\endcases
$$ and $$
b^-(L, p) = \cases
1 & \text{ if } p \rcon L,\cr
-1 & \text{ if } p \lcon L,\cr
0 & \text{ otherwise}.
\endcases
$$  The integers $b^\pm(L, p)$ are well-defined by the definition of full heap.
\enddefinition

The following lemma will be used repeatedly in the sequel, sometimes 
without explicit comment.  

\proclaim{Lemma \sectc.3}
Let $E$ be a full heap over a Dynkin diagram $\Gamma$, let $k$ be a field, 
let $L$ be a finite convex subheap of $E$ and let $p$ be a vertex of $\Gamma$.
Then we have $[H_p, X_L] = c X_L$ and $[H_p, Y_L] = -c X_L$ for some
$c \in \{-2, -1, 0, 1, 2\}$.  More precisely, we have $$\eqalign{
[H_p, X_L] &= b^+(L, p) + b^-(L, p),\cr
[H_p, Y_L] &= - b^+(L, p) - b^-(L, p).\cr
}$$
\endproclaim

\demo{Proof}
This follows from the definition of $H_p$ and a case by case check, similar
to (but easier than) the proof of Lemma \sectb.7.
\qed\enddemo

We define the Weyl group, $W(\Gamma)$, associated to $\Gamma$ to be the
group with generators $\{s_i \in I\}$ indexed by the vertices of $\Gamma$ 
and defining relations $$\eqalign{
s_i^2 &= 1 \text{ for all } i \in I,\cr
s_i s_j &= s_j s_i \text{ if } a_{ij} = 0,\cr
s_i s_j s_i &= s_j s_i s_j \text{ if } a_{ij} < 0 \text{ and } a_{ij} 
a_{ji} = 1,\cr
s_i s_j s_i s_j &= s_j s_i s_j s_i \text{ if } a_{ij} < 0 \text{ and } 
a_{ij} a_{ji} = 2.\cr
}$$  Note that no relation is added in the case where $a_{ij} < 0$ and $a_{ij}
a_{ji} = 4$.

\example{Example \sectc.4}
Define two generalized Cartan matrices $$
A_1 = \left( \matrix
2 & -1 \cr
-2 & 2 \cr
\endmatrix \right) 
\text{\ and\ }
A_2 = \left( \matrix
2 & -2 \cr
-2 & 2 \cr
\endmatrix \right) 
.$$  Then the Weyl group corresponding to $A_1$ is $$
\lan s_1, s_2 : s_1^2 = s_2^2 = 1, (s_1 s_2)^4 = 1 \ran
,$$ isomorphic to the dihedral group of order $8$, and the Weyl group
corresponding to $A_2$ is the infinite group $$
\lan s_1, s_2 : s_1^2 = s_2^2 = 1 \ran
.$$
\endexample

Let $\Pi = \{\a_i : i\in I\}$ and let $\Pi^\vee = 
\{\a_i^\vee : i \in I\}$.  We have a $\zed$-bilinear pairing $\zed\Pi \times
\zed\Pi^\vee \ra \zed$ defined by $$
\lrt{\a_j}{\a_i^\vee} = a_{ij}
,$$ where $(a_{ij})$ is the generalized Cartan matrix.
If $k$ is a field, we extend this to a $k$-bilinear pairing by extension
of scalars.
If $v = \sum_{i \in I} \l_i \a_i$, we write $v \geq 0$ to mean that 
$\l_i \geq 0$
for all $i$, and we write $v > 0$ to mean that $\l_i > 0$ for all $i$.
We view $V = k\Pi$ as the
underlying space of a reflection representation of $W$, determined
by the equalities $s_i(v) = v - \lrt{v}{\a_i^\vee} \a_i$ for all
$i \in I$.

Indecomposable generalized 
Cartan matrices come in three mutually exclusive types.

\proclaim{Theorem \sectc.5 \cite{{\bf 12}, Theorem 4.3}}
Let $A$ be an indecomposable generalized
Cartan matrix.  Then $A$ satisfies one and only 
one of the following three possibilities: 
\item{\rm (i)}{$\det A \ne 0$; there exists $u > 0$ with $Au > 0$; and
$Av \geq 0$ implies $v > 0$ or $v = 0$;}
\item{\rm (ii)}{$A$ has corank $1$; there exists $u > 0$ with $Au = 0$; and
$Av \geq 0$ implies $Av = 0$;}
\item{\rm (iii)}{there exists $u > 0$ with $Au < 0$; and the conditions
$Av \geq 0$ and $v \geq 0$ together imply $v = 0$.}

The matrix $A$ is said to be of finite (respectively, affine, indefinite)
type if it satisfies condition (i) (respectively, (ii), (iii)) above.
\qed\endproclaim

In this paper, we are only concerned with the finite and affine cases above.

Following \cite{{\bf 12}, \S5}, we define a {\it real root} to be a vector
of the form $w(\a_i)$, where $w \in W$ and $\a_i$ is a basis vector.  If
$A$ is of finite type, all roots are real.  If $A$ is of affine type, there
is a unique vector $\d = \sum a_i \a_i$ such that $A\d = 0$ and the $a_i$
are relatively prime positive integers.  Although the notion of {\it imaginary
root} can be defined in general, in the affine type case the imaginary
roots are easily characterized as precisely those vectors of the form 
$n\d$ where $n$ is a nonzero integer.

A {\it root} is by definition a real or imaginary root.  
We denote the set of roots by $\Delta$, as in \cite{{\bf 12}}.
We say a root $\a$ is
positive (respectively, negative) if $\a > 0$ (respectively, $\a < 0$).
If $\a$ is a root, then so is $-\a$, and every root is either positive or
negative.  We will identify the positive (real and
imaginary) roots with elements of $R^+$ as in Definition \sectb.3 so that
$\sum a_i \a_i$ corresponds to the function sending each $i$ to $a_i$.
The {\it height}, $\hei(\a)$ of the root $\a = \sum a_i \a_i$ is by 
definition the integer $\sum a_i$.

\proclaim{Lemma \sectc.6}
Let $A_0$ be a simply laced generalized Cartan matrix of finite type, and let
$\a = \sum a_i \a_i$ and $\be$ be two positive roots associated to $A_0$.  
Define $\a^\vee = \sum a_i \a_i^\vee$, and write $\lrt{\be}{\a^\vee}$ for $
\sum a_i \lrt{\be}{\a_i^\vee}
.$  Then precisely one of the following situations occurs:
\item{\rm (i)}
{$\lrt{\be}{\a^\vee} = 2$ and $\a = \be$;}
\item{\rm (ii)}
{$\lrt{\be}{\a^\vee} = 1$, $\a - \be$ is a root and $\a + \be$ is not a root;}
\item{\rm (iii)}
{$\lrt{\be}{\a^\vee} = 0$ and neither of $\a \pm \be$ is a root;}
\item{\rm (iv)}
{$\lrt{\be}{\a^\vee} = -1$, $\a + \be$ is a root and $\a - \be$ is not a root;}
\item{\rm (v)}
{$\lrt{\be}{\a^\vee} = -2$ and $\a = -\be$.}
\endproclaim

\demo{Proof}
This is well-known, and the proof follows from the argument given in
\cite{{\bf 11}, \S9.4}.
\qed\enddemo

\head \S\sectd.  Parity of heaps in the simply laced case \endhead

Let $A$ be a generalized Cartan matrix,
let $\fg$ be the associated Lie algebra, and
let $\Gamma$ be the corresponding Dynkin diagram.
In \S\sectd, we assume that $A$ is simply laced; in other words,
that $A$ has entries in the set $\{2, 0, -1\}$.  

Let us now fix an orientation of $\Gamma$.

\definition{Definition \sectd.1}
Following \cite{{\bf 12}, (7.8.4)}, we define 
a function $$\sgn : P \times P \ra \{\pm 1\}$$ (depending on the chosen
orientation of $\Gamma$) by the conditions $$
\sgn(p, p') = \cases
-1 & \text{ if } p = p' \text{ or there is an arrow from } p \text{ to } p',\cr
1 & \text{ otherwise.}
\endcases$$  We may extend the above definition to a function
$\sgn : R^+ \times R^+ \ra \zed$ via $$
\sgn(f, g) = \sum_{p \in P} \sum_{q \in P} f(p) g(q) \sgn(p, q)
;$$ similarly, we may extend the definition to a function on
roots by $$
\sgn\left(\sum a_i \a_i, \sum b_j \a_j \right) = 
\sum_i \sum_j a_i b_j \sgn(i, j)
.$$
\enddefinition

\proclaim{Lemma \sectd.2}
Assume additionally that $A$ is of finite type, and that $\a, \be$ are positive
roots such that $\a + \be$ is a root.  Then $\sgn(\a, \be) = -\sgn(\be, \a)$.
\endproclaim

\demo{Proof}
This follows by combining (7.8.7) and (7.8.8) in \cite{{\bf 12}}.
\qed\enddemo

We use the above to define a parity function on finite heaps as follows.

\definition{Definition \sectd.3}
If $F$ is a finite labelled heap over $\Gamma$, we define $$
\sgn(F) = \prod_{{\a, \be \in F} \atop {\a > \be}} \sgn(\e(a), 
\e(\be))
.$$  We extend the notion of parity to finite heaps over $\Gamma$, in the 
obvious way.  (See Example \sectf.4 below for a sample calculation.)
\enddefinition

For our purposes, the key purpose of the $\sgn$-function is the following

\proclaim{Lemma \sectd.4}
Let $\g, \g' \in R^+$, and let $F_1 \in \L_\g$ and $F_2 \in \L_{\g'}$ be 
finite labelled heaps over $\Gamma$.  Then we have $$
\sgn(F_1 \circ F_2) = \sgn(F_1) \sgn(F_2) \sgn(\g', \g)
.$$
\endproclaim

\demo{Proof}
In the computation of the left hand side using Definition \sectd.3, three
types of terms appear in the product: (a) those where $\a$ and $\be$ both
lie in $F_1$, (b) those where $\a$ and $\be$ both lie in $F_2$ and (c)
those where $\a$ lies in $F_1$ and $\be$ lies in $F_2$.  The factorization
in the statement expresses this decomposition.
\qed\enddemo

\definition{Definition \sectd.5}
As in \cite{{\bf 22}}, we define for each $\a \in R^+$ operators $X_\a$ and 
$Y_\a$ on $V_E$ by the formulae $$
X_\a = \sum_{L \in \L_\a(E)} \sgn(L) X_L
$$ and $$
Y_\a = \sum_{L \in \L_\a(E)} \sgn(L) Y_L
.$$
\enddefinition

Although the sums in the above definition may be infinite, note that at
most one summand can act as zero on any given $v_I$.

\proclaim{Lemma \sectd.6}
Maintain the notation of Definition \sectb.6, and let $L \in \L_\a$ and
$L' \in \L_\be$ for some $\a, \be \in R^+$.  If $v_I$ is a basis element
such that $X_L \circ X_{L'}(v_I) \ne 0$, then we have $$
X_L \circ X_{L'}(v_I) = \sgn(\a, \be) X_{L \cup L'}(v_I)
;$$ similarly, if $v_I$ is such that $Y_L \circ Y_{L'}(v_I) \ne 0$ then  $$
Y_L \circ Y_{L'}(v_I) = \sgn(\be, \a) Y_{L \cup L'}(v_I)
.$$
\endproclaim

\demo{Proof}
This is a
consequence of the definitions and Lemma \sectd.4.  The first identity
corresponds to the case $L' \circ L = L \cup L'$, and the second to the 
case $L \circ L' = L \cup L'$.
\qed\enddemo

\head \S\secte.  Representability of roots in the simply laced finite
type case \endhead

In \S\secte, we concentrate on the case of the simply laced, finite type
case.  However, we first require some general results (Definition \secte.1
and lemmas \secte.2\ and \secte.3), which will also be needed in later
sections.

\definition{Definition \secte.1}
If $\a$ is a positive root associated to a Kac--Moody algebra $\fg$, then 
(identifying $\a$ with an element of $R^+$ in the usual way) we call
elements of $\L_\a$ {\it root heaps}.  If $\L_\a$ is 
nonempty, we say that the root $\a$ is {\it representable} in the heap $E$.
\enddefinition

\proclaim{Lemma \secte.2}
Let $\fg$ be a Kac--Moody algebra and let $\a$ be a real non-simple, positive 
root associated to $\fg$.  Then there exists a simple root $\a_i$ such that 
$\lrt{\a}{\a_i^\vee} > 0$ and the root 
$s_i(\a) = \a - \lrt{\a}{\a_i^\vee} \a_i$
is positive.
\endproclaim

\demo{Proof}
This is \cite{{\bf 12}, Proposition 5.1 (e)}.
\qed\enddemo

\proclaim{Lemma \secte.3}
Let $\fg$ be a Kac--Moody algebra with associated Dynkin diagram
$\Gamma$, let $\a = \sum a_i \a_i \in R^+$ be such that $\a > 0$, let
$a = \lrt{\a}{\a_i^\vee}$,
let $E$ be a full heap over $\Gamma$, and let $F \leq E$ with $F \in \L_\a$.
\item{\rm (i)}{If $a = 2$ then $F$ has both a maximal 
vertex $\be$ and a minimal vertex $\be'$ with $\e(\be) = \e(\be') = i$.}
\item{\rm (ii)}{If $a = 1$ then either $F$ has a maximal 
vertex $\be$ with $\e(\be) = i$, or minimal vertex $\be'$ with 
$\e(\be') = i$, but not both.}
\item{\rm (iii)}{If $a = 0$ and $F$ has a maximal element labelled $i$, then
there exists $\be \in E \backslash F$ with $\e(\be) = i$ such that 
$F' = F \cup \{\be\}$ is convex and $\be$ is minimal in $F'$.}
\item{\rm (iv)}{If $a = 0$ and $F$ has a minimal element labelled $i$, then
there exists $\be \in E \backslash F$ with $\e(\be) = i$ such that 
$F' = F \cup \{\be\}$ is convex and $\be$ is maximal in $F'$.}
\endproclaim

\demo{Proof}
Let us write $\a = \sum a_k \a_k$ as a sum of simple roots.  Since
$a > 0$, we must have $a_i > 0$ for some $i$; in other words, $F$ contains 
at least one
element labelled $i$.  Let $\z_0$ and $\z_1$ denote the least and greatest
elements of $F \cap \e^{-1}(i)$, respectively.  

In order to calculate $\lrt{\a}{\a_i^\vee}$, the only relevant summands in
the expression for $\a$ are 
those corresponding to $\a_i$ itself and to the simple roots adjacent to 
$\a_i$.  Define $F'$ to be the set of all $\g \in F$ with $\e(\g)$ adjacent
to $i$.
The definition of full heap shows that there are three possibilities
for elements $\g \in F'$: (a) $\g < \z_0$, (b) $\g > \z_1$,
or (c) $\g$ lies between two elements $\z$, $\z'$ of $F$ with 
$\e(\z) = \e(\z') = i$ and $(\z, \z') \cap \e^{-1}(i) = \emptyset$.

Let us first consider case (c).
If $a_{ij} = -2$, such an open interval $(\z, \z')$ contains a unique element
$\g$ with label adjacent to $i$.  The other possibility is that $a_{ij} = -1$, 
in which case $(\z, \z')$ contains precisely two elements, $\g$ and $\g'$, 
with labels adjacent to $i$.  Let $\a'$ be the element of $R^+$ given by the
character $\cha([\z_0, \z_1])$ of the closed interval $[\z_0, \z_1]$.  The 
above case analysis in terms of $a_{ij}$ shows that $\lrt{\a'}{\a_i^\vee} = 
\lrt{\a_i}{\a_i^\vee} = 2$, and this identity also holds in the case that 
$[\z_0, \z_1]$ consists of a single element.

The contributions to $\lrt{\a}{\a_i^\vee}$ that do not come from
$\lrt{\a'}{\a_i^\vee}$ must therefore come from the elements $\g$ in cases
(a) and (b) above, which means in particular that $\lrt{\a}{\a_i^\vee}
\leq \lrt{\a'}{\a_i^\vee}$.  

If $a = 2$, then there cannot be any such elements, or the latter 
inequality would be strict, contrary to hypothesis.  This
means that $\z_0$ (respectively, $\z_1$) is minimal (respectively, maximal)
in $F$, which establishes assertion (i).

If $a = 1$, then there must be precisely one such element $\g$ (and, in fact,
$a_{ij} = -1$ must also hold).  This means that either $\z_0$ is minimal
in $F$, or $\z_1$ is maximal in $F$, but not both, which establishes 
assertion (ii).

Suppose $a = 0$ and that there is a maximal element labelled $i$.  There are
then no elements $\g$ corresponding to case (a).  Furthermore, either there
is precisely one element $\g$ arising from case (b) and $a_{ij} = -2$, or
there are precisely two elements $\g$ arising from case (b) and $a_{ij} = -1$.
If we write $\z_0 = E(i, t)$, then setting $\be = E(i, t-1)$ will then 
satisfy the hypothesis of (iii) by the definition of full heap.

Part (iv) is proved by a symmetrical argument.
\qed\enddemo

We now let $A$ be a simply laced generalized Cartan matrix of (necessarily
untwisted) affine type, corresponding to the finite type matrix $A_0$.
The next result, whose method of proof is familiar from Kashiwara's
celebrated Grand Loop \cite{{\bf 13}, \S4}, establishes the basic properties of 
root heaps and the operators of Definition \sectd.5.

\proclaim{Proposition \secte.4}
Let $E$ be a full heap over the Dynkin diagram $\Gamma$ of $A$, 
let $A_0$ be the corresponding finite type generalized Cartan matrix 
with Kac--Moody algebra $\fg_0$, and let $\a = \sum \l_i \a_i$ be a positive 
real root associated to $\fg_0$.
\item{\rm (i)}{The root $\a$ is representable in $E$.}
\item{\rm (ii)}{The operator $X_\a$ is nonzero, lies in the Lie algebra
generated by the $X_p$, and (in the case $k = \complex$) is equal to the 
element $E_\a$ in the notation
of \cite{{\bf 12}, (7.8.5)}.}
\item{\rm (iii)}{If $p$ is any vertex of $\Gamma$, then 
$[H_p, X_\a] = \lrt{\a}{\a_p^\vee} X_\a$.}
\item{\rm (iv)}{The operator $Y_\a$ is nonzero, lies in the Lie algebra
generated by the $Y_p$, and (in the case $k = \complex$) is equal to the 
element $-E_{-\a}$ in the
notation of \cite{{\bf 12}, (7.8.5)}.}
\item{\rm (v)}{If $p$ is any vertex of $\Gamma$, then 
$[H_p, Y_\a] = - \lrt{\a}{\a_p^\vee} Y_\a$.}
\item{\rm (vi)}{For any proper ideal $I$ of $E$, there do not exist root
heaps $L, L' \in \L_\a$ such that both $I \cup L$ and $I \backslash L'$ are
ideals.}
\item{\rm (vii)}{We have $H_\a = \sum \l_i H_{\a_i} = \a^\vee$, where $\a^\vee$
is as defined in \cite{{\bf 12}, \S5.1}.}
\endproclaim

\demo{Proof}
We will prove the statements simultaneously by induction on $\hei(\a)$.
The proofs of (iv) and (v) are very similar to those of (ii) and (iii),
respectively, so we do not include them.

The base case is $\hei(\a) = 1$, in other words, $\a$ is simple.  Parts (i)
and (ii) follow from part (a) of the definition of a fibred heap and the
definitions of \cite{{\bf 12}, \S7.8}, part
(iii) is immediate from Theorem \sectc.1, part (vi) follows from
the definition of a full heap and part (vii) is trivial (again using
the definitions of \cite{{\bf 12}, \S7.8}).

For the inductive step, we use Lemma \secte.2 to find a simple root $\a_i$ 
such that $a = \lrt{\a}{\a_i^\vee} > 0$ and the root 
$\a' = s_i(\a) = \a - a \a_i$ is positive.  By the
inductive hypothesis, we have $$
[H_i, X_{\a'}] = \lrt{\a'}{\a_i^\vee}
= \lrt{\a - \lrt{\a}{\a_i^\vee} \a_i}{\a_i^\vee} = 
\lrt{\a}{\a_i^\vee}(1 - \lrt{\a_i}{\a_i^\vee})
= - \lrt{\a}{\a_i^\vee}
.$$  Since $\a'$ is representable by the inductive hypothesis, 
Lemma \sectc.3 shows that we have $\lrt{\a}{\a_i^\vee} \in \{1, 2\}$.  

This gives three possibilities for a root heap $L \in \L_{\a'}$.  
Case (a) is that $a = 2$, $i \lcon L$ and $L \rcon i$, but in fact this
cannot occur, because it implies (by the inductive hypothesis) that 
$\lrt{\a'}{\a_i^\vee} = -2$, which contradicts Lemma \sectc.6.
Case (b) is that $a = 1$, and $i \lcon L$ and but we do not
have either $L \rcon i$, $L \lcon i$ or $i \rcon L$.
Case (c) is that $a = 1$, and $L \rcon i$
but we do not have either $i \lcon L$, $i \rcon L$ or $L \lcon i$.

We are now in case (b) or (c), so that $a = 1$.  
In this situation, $\a$ is representable because we can add a new 
maximal or minimal vertex labelled $i$ to $L$ to form a root
heap in $\L_\a$.  By Lemma \sectb.1 (vii), this means that 
$X_\a$ is nonzero, and by Lemma 
\secte.3, any root heap in $\L_{\a}$ has 
either a maximal or minimal vertex (but not both) labelled $i$.
We claim in this case that $X_\a = \pm [X_{\a_i}, X_{\a'}]$.  It is enough
to check that each side of the equation acts in the same way on a basis
vector, $v_I$.  The right hand side of the equation is equal to $$
X_{\a_i} \circ X_{\a'} 
- X_{\a'} \circ X_{\a_i}
.$$    It now follows that unless we have $L' \succ I$ for some $L' \in \L_\a$,
both sides of the equation will act as zero, so let us assume that this
condition is satisfied.  By Lemma \secte.3 (ii), every element 
$L'$ of $\L_\a$ is
uniquely of the form $\{\be\} \circ L$ or of the form $L \circ \{\be\}$ (but
not both) for some $\be$ and $L \in \L_{\a'}$, and we have shown that any 
such $L \in \L_{\a'}$ can be extended to an element of $\L_\a$ in this way.

In case (b), $X_{\a_i} \circ X_{\a'}$ acts as zero on $v_I$, so by Lemma
\sectd.6 we have $$
[X_{\a_i}, X_{\a'}] . v_I = -X_{\a'} \circ X_{\a_i} v_I = - \sgn(\a', \a_i)
X_\a . v_I
.$$
In case (c), $X_{\a'} \circ X_{\a_i}$ acts as zero on $v_I$, so by Lemma
\sectd.6 we have $$
[X_{\a_i}, X_{\a'}] . v_I = X_{\a_i} \circ X_{\a'} v_I = \sgn(\a_i, \a')
X_\a . v_I
.$$  Lemma \sectd.2 now shows that $[X_{\a_i}, X_{\a'}] = \sgn(\a_i, \a')
X_{\a}$.  Since, by \cite{{\bf 12}, (7.8.5)}, we have $[E_{\a_i}, E_{\a'}] 
= \sgn(\a_i, \a') E_{\a}$ (a formula also valid for negative roots), we
have $X_\a = E_\a$ by the inductive hypothesis, completing the proof of
(i) and (ii).

We now observe, using the Jacobi identity and the inductive hypothesis,
that $$\eqalign{
[H_p, X_\a] 
&= \sgn(\a', \a_i) [H_p,  [X_{\a_i}, X_{\a'}]] \cr
&= \sgn(\a', \a_i) \left( [H_p, X_{\a_i}], X_{\a'}] 
+ [X_{\a_i}, [H_p, X_{\a'}]] \right) \cr
&= \sgn(\a', \a_i) (\lrt{\a_i}{\a_p^\vee} + \lrt{\a'}{\a_p^\vee})
[X_{\a_i}, X_{\a'}] \cr
&= \lrt{\a}{\a_p^\vee} \sgn(\a', \a_i) [X_{\a_i}, X_{\a'}] \cr
&= \lrt{\a}{\a_p^\vee} X_\a,
}$$ which completes the proof of (iii).

To prove (vi), we write $\a = \a' + \a_i$ as before, so that 
$\lrt{\a}{\a_i^\vee} = 1$ by Lemma \sectc.6.  Let $I, L, L' \in \L_\a$ be as
in the statement, contrary to hypothesis.  By the inductive hypothesis 
applied to (iii), we have $[H_i, X_{\a'}] = X_{\a'}$.  This gives two 
possibilities for $L$ (and similar possibilities for $L'$): 
the first is that $i \rcon L$ and but we do not
have either $L \lcon i$, $L \rcon i$ or $i \lcon L$, and the second is that 
$L \lcon i$ but we do not have either $i \rcon L$, $i \lcon L$ or $L \rcon i$.
Let us consider the first possibility.  Because
$\cha(L) = \cha(L')$, $L'$ must contain an element labelled $i$, and since
$L \cup L'$ is convex (as $I \cup L$ and $I \backslash L'$ are ideals),
we must have $L' \rcon i$.  This is not a permissible
configuration for $L'$, so we have a contradiction and we conclude that
in fact $L \lcon i$.  A dual analysis shows that $i \rcon L'$.  If we
now delete the maximal element in $L$ with label $i$ to form a heap $L_0$,
and we delete the minimal element in $L'$ with label $i$ to form a heap $L'_0$,
then the inductive hypothesis applied to the ideal $I$ and the heaps
$L_0, L'_0 \in \L_{\a'}$ shows that this situation is impossible,
proving (vi).

It follows from (vi) that $[X_\a, Y_\a] = H_\a$,
and we know from \cite{{\bf 12}, \S7.8} that $[E_\a, E_{-\a}] = -\a^\vee$,
so we have $H_\a = \a^\vee$ by (ii) and (iv).  The other assertions follow
from \cite{{\bf 12}, (5.1.1), \S7.8}, using the fact that all roots have the
same length.  This establishes (vii).
\qed\enddemo

\proclaim{Corollary \secte.5}
Let $E$ be a full heap over the Dynkin diagram $\Gamma$ of $A$, 
and let $A_0$ be the corresponding finite type generalized Cartan matrix.
If $\a, \be, \g$ are positive roots associated to $A_0$ such that
$\a = \be + \g$, then any root heap $L \in \L_\a$ decomposes uniquely as
a disjoint union $L = L_1 \cup L_2$ of (convex) subheaps $L_1 \in \L_\be$
and $L_2 \in \L_\g$ such that one of $L_1$ and $L_2$ is an ideal of $L$ and
the other is a filter of $L$.
\endproclaim

\demo{Proof}
We can choose a proper ideal $I$ of $E$ such that $L \prec I$ by Lemma
\sectb.1 (vii).
We have $H_\a (v_I) = v_I$, and by Proposition \secte.4 (vii),
we have $H_\a = H_\be + H_\g$.  By Proposition \secte.4 (vi), there are 
two ways this can happen: either $H_\be(v_I) = v_I$ and $H_\g(v_I) = 0$,
or vice versa.  In the first case, there is a filter $L_1$ of $L$ and an ideal
$L_2$ of $L$ with $L_1 \in \L_\be$ and $L_2 \in \L_\g$, and in the second case,
there is an ideal $L_1$ of $L$ and a filter $L_2$ of $L$ with $L_1 \in \L_\be$
and $L_2 \in \L_\g$.
\qed\enddemo

\head \S\sectf. The non simply laced case \endhead

The methods presented for simply laced Lie algebras can be generalized
to the non simply laced case.  The right way to consider this seems to be
to regard the non simply laced objects as folded versions of their simply
laced untwisted affine counterparts.  For roots, this is the procedure 
described in \cite{{\bf 12}, \S7.9}.  For our purposes, we also need a
categorified version of this phenomenon suitable for full heaps.

Let $A$ be a simply laced generalized Cartan matrix of untwisted affine type 
and let $\Gamma$ be the corresponding Dynkin diagram, and suppose that $\mu$
is a nonidentity graph automorphism of $\Gamma$.  Although in general
$\mu$ can be of order $2$ or $3$, we will make two additional assumptions about
$\mu$: (a) $\mu$ has order precisely $2$ and (b) for any vertex $p$, $\mu(p)$ 
and $p$ are not distinct adjacent vertices.

The group $\{1, \mu\}$ acts on the Dynkin
diagram $\Gamma$, and we denote the orbit containing the vertex $p$ by
$f(p) = \bar{p}$.  This induces an action on the simple roots $\a_i$, and we
extend this to a linear action on $k \otimes_\zed R$ by $\mu(a_i \a_i) = 
a_i \mu(\a_i)$.
Let us also define $f(\a_i) = (\a_i + \mu(\a_i))/2$ and extend linearly
to $k \otimes_\zed R$.

Let $A$ and $\Gamma$ be as above and $\Delta$ be the set of roots for $A$.
It is known \cite{{\bf 12}, Proposition 7.9} that 
the set $\{f(\a) : \a \in \Delta\}$ is a root system for a Kac--Moody
algebra $\overline{\fg}$ with 
simple roots $\{f(\a_i)\}$.  A root $f(\a)$ is called {\it long}
if $\a = \mu(\a)$, and {\it short} otherwise.
The Dynkin diagram $\overline{\Gamma}$ for $\overline{\fg}$ has vertices
labelled by the orbits $\bar{p}$, and is such that if $p$ and $q$ are
distinct vertices of $\Gamma$, then $p$ and $q$
are adjacent in $\Gamma$ if and only if the (distinct) vertices $\bar{p}$
and $\bar{q}$ are adjacent in $\overline{\Gamma}$.  If $\Gamma$ contains
three vertices $p$, $\mu(p)$ and $q$ such that $q$ is adjacent to both
$p$ and $\mu(p)$, then we join $\bar{p}$ and $\bar{q}$ in $\overline{\Gamma}$
by a double edge with
an arrow pointing towards $\bar{p}$.  (It is possible for this procedure to
result in a double edge with two arrows in opposite directions.)

We will say that $A$ (respectively, 
$\Gamma$) {\it folds} to $\overline{A}$ (respectively, $\overline{\Gamma}$) 
via $\mu$.

\proclaim{Proposition \sectf.1}
Let $E = (E, \leq, \e)$ be a full (labelled) heap over the Dynkin diagram 
$\Gamma$ of $A$, where $A$ is a simply laced generalized Cartan matrix of 
untwisted affine type.  
Suppose that $\mu$ folds $A$ and $\Gamma$ to $\overline{A}$ and
$\overline{\Gamma}$, and also that whenever we have vertices $p, q$
of $\Gamma$ satisfying (a) $\mu(p) \ \C \  q$, (b) $\a \in \e^{-1}(p)$ and 
(c) $\be \in \e^{-1}(q)$, then $\a$ and $\be$ are comparable in $E$.
Then $\overline{E} = (E, \leq, f \circ \e)$
is a full (labelled) heap over $\overline{\Gamma}$.
\endproclaim

\demo{Proof}
We first show that $\overline{E}$ is a heap, \idest that Definition
\secta.1 holds.  

For part 1 of the definition, 
it is enough to show that if $\a, \be \in E$, then
$\a$, $\be$ are comparable if $f \circ \e(\a) = f \circ \e(\be)$.
There are four cases to consider: either $\e(\a) = \e(\be)$, or $\e(\a)
= \mu(\e(\be))$, or $\e(\a)$ is adjacent to $\e(\be)$, or $\e(\a)$ is
adjacent to $\mu(\e(\be))$.  In each case, $\a$ and $\be$ are guaranteed
to be comparable either by the definition of a heap, or by the hypotheses
on $\mu$ in the statement.
 
For part 2, suppose that $\a, \be \in E$, $\a \leq \be$ and $\e(\a) \ \C \ 
\e(\be)$.  It is immediate that $f(\e(\a)) \ \C \ f(\e(\be))$, from which
the assertion follows.

We next show that $\overline{E}$ is fibred.  Part (a) of Definition
\secta.6 comes from the fact that $(f \circ \e)^{-1}(\bar{p}) \subseteq
\e^{-1}(p)$.  For (b), assume that $p'$ and $q'$ are adjacent vertices of
$\Gamma$ and that $\a \in E$ satisfies $f(\e(\a)) = p'$.  The properties
of graph automorphisms guarantee the existence of a vertex $q$ of $\Gamma$
adjacent to $p = \e(\a)$ with $f(q) = q'$.  Since $E$ is fibred, there
exists $\be \in E$ such that $\a$ covers $\be$ or $\be$ covers $\a$, with
$\e(\be) = q$.  This vertex satisfies $f(\e(\be)) = q'$, as required.

Finally, we show that $\overline{E}$ is full.  
Let $\a, \be \in E$ (where $\a < \be$)
be such that $f \circ \e(\a) = f \circ \e(\be) = \bar{p}$ and 
$(\a, \be) \cap (f \circ \e)^{-1}(\bar{p}) = \emptyset$.  

Suppose first that
$\e(\a) = \e(\be) = p,$ say.  In this case, the interval $(\a, \be)$ in $E$
contains precisely two elements, $\g, \g'$, with labels adjacent to $p$.  
This shows that the elements $\g$ and $\g'$, considered as elements of 
$\overline{E}$, are the only two elements of $(\a, \be)$ with labels
adjacent to $\bar{p}$.  It remains to show that the Dynkin diagram
$\overline{\Gamma}$ does not contain an arrow from $f \circ \e(\g)$ or
$f \circ \e(\g')$ to $\bar{p}$; we deal with the former case, the other being
similar.  
Now either $\be$ covers $\g$, or $\g$ covers $\a$ (or possibly both); 
we prove the former, and the latter follows by a dual argument.
If such an arrow exists, there must exist
$\be'$ with $\e(\be') = \mu(\e(\be))$ and either $\be'$ covers $\g$ or
$\g$ covers $\be$.  If $\be'$ covers $\g$, this implies that $\be$ and $\be'$
are comparable, contrary to the hypothesis on $\mu$.  On the other hand, if
$\g$ covers $\be'$ then by hypothesis, $\be'$ and $\a$ are comparable, which
forces $\a < \be' < \be$.  This in turn implies that $$
\be' \in (\a, \be) \cap (f \circ \e^{-1}) (\bar{p})
,$$ a contradiction.

The other case to consider is that $\e(\a) \ne \e(\be)$, which implies that 
$\mu(\e(\a)) = \e(\be)$.  If $\a = E(p, t)$ in the numbering of Remark
\secta.7, define $\a' = E(p, t+1) > \a$.  By the hypothesis on $\mu$,
$\a'$ and $\be$ are comparable, which (by the hypotheses on $\a$ and $\be$) 
means that $\a < \be < \a'$ and that $(\be, \a')$ is nonempty.  Since the
interval $(\a, \a')$ in $E$ contains precisely two elements, $\g$ and $\g'$,
with labels adjacent to $\e(\a)$, and neither label is $\e(\be)$, we may
assume without loss of generality that $\g \in (\a, \be)$ and $\g' \in
(\be, \a')$.  This means that the interval $(\a, \be)$ in $\overline{E}$
contains precisely one vertex, $\g$, with label adjacent to $\bar{p}$.
It remains to check that there is an arrow in the Dynkin diagram
from $f \circ \e(\g)$ to $\bar{p}$, and this follows from that fact that
$\be$ covers $\g$ and $\g$ covers $\a$.
\qed\enddemo

\remark{Remark \sectf.2}
The words ``labelled'' may be dropped from the statement of Proposition
\sectf.1 using a familiar argument.
In the situation of the proposition, we will say that {\it $E$ folds to 
$\overline{E}$}.
\endremark

All the examples we know of full heaps over non simply laced Dynkin diagrams
for affine Kac--Moody algebras
are obtained from the simply laced Dynkin diagrams by the folding procedure
just described.  (See the Appendix for details.)

\definition{Definition \sectf.3}
Suppose that $E = (E, \leq, \e)$ is a full heap over the Dynkin diagram 
$\Gamma$ of $A$, where $A$ is a simply laced generalized Cartan matrix of 
untwisted affine type, and that $\mu$ is a diagram automorphism of $\Gamma$
that folds the triple $(A, \Gamma, E)$ to $(\overline{A}, \overline{\Gamma},
\overline{E})$.  An orientation of $\Gamma$ is said to be {\it compatible}
with $\mu$ if $\sgn(p, p') = \sgn(\mu(p), \mu(p'))$.  
If $L$ is a finite subheap of $\overline{E}$, and $\Gamma$ has an 
orientation compatible with $\mu$, then we define $\overline{L}$ to be the
subheap of $\overline{E}$ corresponding to $L$, and we define
$\sgn(\overline{L}) = \sgn(L)$, where parity is taken with respect to 
this compatible orientation.  (Every finite subheap of $\overline{E}$ arises
in this way.)

The operators $X_\a$ and $Y_\a$ in the non simply laced case are now defined
in the same way as in Definition \sectd.5.  (The arrows induced by the 
orientation have nothing to do with the arrows used in the definition of 
the Dynkin diagram.)
\enddefinition

\example{Example \sectf.4}
Let $\Gamma$ be the Dynkin diagram of type $A_5^{(1)}$, and let $E$ be the
full heap over $\Gamma$ shown in Figure 4 of the Appendix.  
Let $\overline{E}$ be the 
corresponding heap over the Dynkin diagram $\overline{\Gamma}$ of type
$C_3^{(1)}$ shown in Figure 6.  Suppose the Dynkin diagrams are oriented
as in Figure 1, and let $F$ be a convex subheap of $E$ with character
$\a_2 + \a_3 + \a_4 + \a_5$ (all such subheaps are isomorphic).
Then the subheap $\overline{F}$ of $\overline{E}$ corresponding to $F$
has character $\a_1 + 2\a_2 + \a_3$.  (See Figure 2.)  There are two
pairs $(\a, \be)$ of elements in $F$ such that $\a > \be$ and either 
$\sgn(\e(\a), \e(\be)) = -1$ or $\e(\a) = \e(\be)$.  These are the pairs 
$(\a, \be)$ and $(\g, \be)$, where $\e(\a) = 5$, $\e(\be) = 4$ and 
$\e(\g) = 3$.  We thus have $\sgn(F) = (-1)^2 = 1$, and we have 
$\sgn(\overline{F}) = \sgn(F)$ by definition.
\endexample

\topcaption{Figure 1} Compatible orientations for the Dynkin diagrams of types
$A_5^{(1)}$ and $C_3^{(1)}$ in Example \sectf.4
\endcaption 
\centerline{
\hbox to 4.638in{
\vbox to 1.125in{\vfill
        \includegraphics{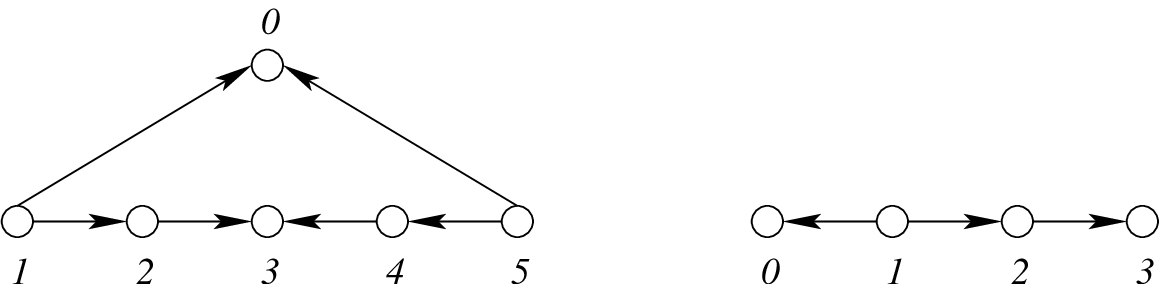}
}
\hfill}
}

\topcaption{Figure 2} The heaps $F$ and $\overline{F}$ of Example 
\sectf.4
\endcaption 
\centerline{
\hbox to 2.222in{
\vbox to 0.902in{\vfill
        \includegraphics{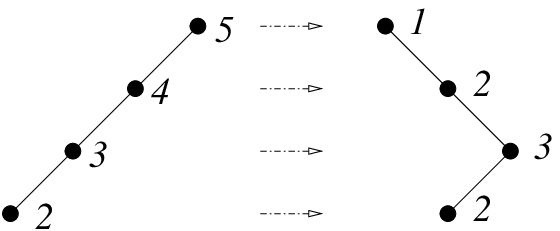}
}
\hfill}
}

It is immediate from the definitions that if $E$ folds to $\overline{E}$
via $\mu$ and $L \in \L_\a$ is a root heap in $E$, then
we have $\overline{L} \in
\L_{f(\a)}$.  The following result shows that the converse is also true, so
that we may pass easily between the root heaps of $E$ and those of 
$\overline{E}$.

\proclaim{Proposition \sectf.5}
Let $E$ be a full heap over a simply laced Dynkin diagram $\Gamma$ of 
untwisted affine type, $A$, 
let $A_0$ be the corresponding finite type generalized Cartan matrix 
with Kac--Moody algebra $\fg_0$, and let $\a = \sum \l_i \a_i$ be a positive 
real root associated to $\fg_0$.  
Suppose that the map $f$ sends the roots of $\fg_0$ to
the roots of another Kac--Moody algebra $\overline{\fg_0}$ of finite type, 
identifying simple roots with simple roots,
and that $E$ folds to $\overline{E}$
via $\mu$.  Let $f(\a)$ be the root of the simple Lie algebra
$\overline{\fg_0}$ corresponding to $\a$, and
assume that the field $k$ does not have characteristic $2$.
\item{\rm (i)}{The root $f(\a)$ is representable in $\overline{E}$, and 
for any root heap $\overline{L} \in \L_{f(\a)}$, we have
$L \in \L_\a \cup \L_{\mu(\a)}$.}
\item{\rm (ii)}{The operator $X_\a$ is nonzero, lies in the Lie algebra
generated by the operators $X_p$ on $V_{\overline{E}}$, 
and (in the case $k = \complex$) is equal to the element $E_\a$ in the notation
of \cite{{\bf 12}, (7.9.3)}.}
\item{\rm (iii)}{The operator $Y_\a$ is nonzero, lies in the Lie algebra
generated by the operators $Y_p$ on $V_{\overline{E}}$, 
and (in the case $k = \complex$) is equal to the 
element $-E_{-\a}$ in the
notation of \cite{{\bf 12}, (7.9.3)}.}
\endproclaim

\demo{Proof}
As in Proposition \secte.4, we prove the statements by simultaneous
induction on $\hei(f(\a))$, calculated with respect to the basis of simple
roots $f(\a_i)$.  The case where $\a$ is simple follows from the definitions.

Suppose now that $f(\a) \ne f(\a_i)$ for any simple root $\a_i$.  It follows
from the definitions that the set $\L_{f(\a)}$ is nonempty: we may choose 
$L_1 \in \L_\a$ by Proposition \secte.4 (i), and then $\overline{L_1} \in
\L_{f(\a)}$.  

Now let $\overline{L} \in \L_{f(\a)}$ be arbitrary.  Since
$\overline{E}$ is full by Proposition \sectf.1, we may apply Lemma
\secte.2 to find a simple root $f(\a_i)$ such that $\lrt{f(\a)}{f(\a_i)}
> 0$; we may then apply Lemma \secte.3 to $\overline{E}$.  Let $I$
be a proper ideal of $E$ with $L \succ I$; this exists by Lemma \sectb.1
(vii).  

Suppose that we are in case (i) of Lemma \secte.3, so that 
$\overline{L}$ has a maximal vertex
$\be$ and a minimal vertex $\be'$ with $f \circ \e(\be) = f \circ \e(\be')
= \overline{i}$.  
Since we are in case (i), we have $$
s_{\overline{i}}(f(\a)) = f(\a) - 2 f(\a_i) = f(\a) - \a_i - \mu(\a_i) = f(\a')
,$$ for some root $\a'$.
This means that $\overline{L} \backslash \{\be, \be'\} \in \L_{f(\a')}$,
which by the inductive hypothesis shows that $L \backslash \{\be, \be'\}
\in \L_{\a'}$.  
We cannot have $\e(\be) = \e(\be')$, as this would contradict Lemma \sectc.3,
Lemma \sectc.6 and Proposition \secte.4 (iii) applied to $L \backslash
\{\be, \be'\}$.
It must therefore be the case that $p = \e(\be) = \mu(\e(\be')) = q$, where 
$\mu(\e(\be))$ and
$\e(\be)$ are distinct.  
Now consider the element of $V_E$ given by $$
[X_p, [X_q, X_{\a'}]] . v_I
.$$  The bracketed expression expands to $$
X_p \circ X_q \circ X_{\a'}
- X_p \circ X_{\a'} \circ X_q
- X_q \circ X_{\a'} \circ X_p
+ X_{\a'} \circ X_q \circ X_p
.$$  Since $\mu$ is compatible with $E$, any elements of
$E$ with labels $p$ and $q$ are comparable, but since $p$ and $q$ are not
adjacent, we have $X_p \circ X_q = X_q \circ X_p = 0$.  Since $L$ has no
minimal element labelled $p$, we must have $X_p . v_I = 0$.  It follows 
that $$
[X_p, [X_q, X_{\a'}]] . v_I = - X_p \circ X_{\a'} \circ X_q . v_I
,$$ which is a nonzero multiple of $v_{I \cup L}$.  By Proposition 
\secte.4 (ii) and the properties of the Chevalley bases given in \cite{{\bf 12},
(7.8.5)}, we see that $[E_p, [E_q, E_{\a'}]]$ must be a nonzero multiple
of $E_{\a' + \a_i + \mu(\a_i)}$, so that in particular $\a' + \a_i + 
\mu(\a_i)$ is a root and $L \in \L_\a$, proving (i).  (A similar argument shows
that $\a' + \a_i$ and $\a' + \mu(\a_i)$ are both roots.)  To prove (ii), we 
note that $f(\a_i)$ is a short root; furthermore, because $\a' + \a_i$ and 
$\a' + \mu(\a_i)$ are roots, $\mu(\a') + \a_i$ must also be a root.  By
\cite{{\bf 12}, (7.9.6)}, if $\g$ and $\g'$ are short roots whose sum is a root,
then $\g + \mu(\g')$ is not a root.  It follows that $\a'$ must be a long root,
which by the choice of orientation on $\Gamma$ means that $$
\sgn(\a_i, \a' + \mu(\a_i)) = 
\sgn(\mu(\a_i), \a' + \a_i)
.$$  Applying Lemma \sectd.6 now shows that $\sgn(L) = -\sgn(L \backslash
\{\be, \be'\})$.  As operators on $V_{\overline{E}}$, we have $$\eqalign{
[X_{\bar{p}}, [X_{\bar{p}}, X_{f(\a')}]] . v_I
&= -2 X_{\bar{p}} \circ X_{f(\a')} \circ X_{\bar{p}} . v_I \cr
&= 2 X_{f(\a)} . v_I, \cr
}$$ where we have used the fact that $X_{\bar{p}} \circ X_{\bar{p}}$ is zero.
Since every element of $\L_{f(\a)}$ has a maximal and a minimal element 
labelled $\bar{p}$, we see that $$
[X_{\bar{p}}, [X_{\bar{p}}, X_{f(\a')}]] = 2 X_{f(\a)}
.$$  A similar calculation using the Chevalley basis \cite{{\bf 12}, (7.9.3)} shows
that $$
[E_{\bar{p}}, [E_{\bar{p}}, E_{f(\a')}]] = 2 \sgn(\a_i, \a') 
\sgn(\mu(\a_i), \a' + \a_i) E_{f(\a)} = 2 E_{f(\a)}
,$$ thus proving (ii).

The other possibility is that we are in case 
(ii) of Lemma \secte.3, so that $$
s_{\overline{i}}(f(\a)) = f(\a) - f(\a_i) = f(\a')
.$$  We will deal with the subcase where $\overline{L}$ has a minimal vertex
$\be$ labelled $\bar{p}$, but no such maximal vertex.
Let us assume that $\e(\be) = p$.  A similar, but
easier, argument establishes that $\overline{L} \backslash \{\be\} \in
\L_{f(\a')}$ and $L \backslash \{\be\} \in \L_{\a'}$.  
If $\overline{L}$ has a maximal vertex labelled $\bar{p}$ but no such minimal
vertex, a similar argument holds.
Operating on $V_E$, we then find (by acting both sides on a suitable $v_I$)
that $$
[X_p, X_{\a'}] = \sgn(\a_p, \a') X_{\a}
$$ by Lemma \sectd.6, Proposition \secte.4 (ii) and \cite{{\bf 12}, (7.8.5)},
from which (i) follows. 
Analogous calculations on $V_{\overline{E}}$ then
show that $$
[X_{\bar{p}}, X_{f(\a')}] = \sgn(\a_p, \a') X_{f(\a)}
,$$ proving (ii).

The proof of (iii) follows by symmetric arguments.  In the first case
above, this results in the identities $$
[Y_{\bar{p}}, [Y_{\bar{p}}, Y_{f(\a')}]] = 2 Y_{f(\a)}
$$ and $$
[E_{\bar{p}}, [E_{\bar{p}}, E_{f(\a')}]] = 2 E_{f(\a)}
.$$  In the other case, we obtain $$
[Y_{\bar{p}}, Y_{f(\a')}] = \sgn(\a', \a_p) Y_{f(\a)} = -\sgn(\a_p, \a') 
Y_{f(\a)}
.$$
\qed\enddemo

\definition{Definition \sectf.6}
Let $A$ be a generalized Cartan matrix of untwisted affine type
with Dynkin diagram $\Gamma$.
If either
\item{(i)}{$A$ is simply laced and $E$ is any full heap over 
$\Gamma$, or }
\item{(ii)}{$A$ is not simply laced and occurs as
a matrix $\overline{A}$ arising from a folded heap $E = \overline{E'}$ as in
Proposition \sectf.1,}

\noindent then we call $E$ a {\it simply folded} full heap over $\Gamma$.
\enddefinition

When restricted to the simply laced case, the following theorem is similar 
to the unproven \cite{{\bf 22}, Theorem 4.1}.

\proclaim{Theorem \sectf.7}
Let $E$ be a simply folded full heap over the Dynkin diagram $\Gamma$ 
of the generalized
Cartan matrix $A$ of an untwisted affine Kac--Moody algebra.

Let $A_0$ be the corresponding finite type generalized Cartan matrix 
with Kac--Moody algebra $\fg_0$ and set of positive roots $\Delta^+$.  Then 
the set of operators $$
\{X_\a : \a \in \Delta^+\} \cup
\{Y_\a : \a \in \Delta^+\} \cup
\{H_p: p \text{ is a vertex of } \Gamma\}
$$ on $V_E$ over the field $k = \complex$ is linearly independent and its 
span is isomorphic to the (simple) Lie algebra $\fg_0$; in particular, the
isomorphism type depends only on $\fg_0$ (rather than $E$).
\endproclaim

\demo{Proof}
From \cite{{\bf 12}, \S7.8}, we know that, over $k = \complex$, the algebra 
$\fg_0$ has a basis given
by $$\{E_\a : \a \in \Delta\} \cup \{E_{-\a} : \a \in \Delta\} \cup
\{ \a_p^\vee : p \text{ is a vertex of } \Gamma\}.$$  The conclusion now 
follows from Theorem \sectc.1, Proposition \secte.4 (ii) and (iv) (in
the simply laced case) and Proposition \sectf.5 (ii) and (iii) (in the non
simply laced case).
\qed\enddemo

This theorem can be used to construct all finite dimensional simple Lie 
algebras over $\complex$ except those of types $E_8$, $F_4$ and $G_2$.
As we explain in \S\secti, for the simple Lie algebras other than these three, 
it is possible
to perform the construction using a finite dimensional subspace of $V_E$, and
this leads to combinatorial constructions of the spin modules in types $B$ 
and $D$ without using Clifford algebras.  (The type $D$ construction has
already been described without proof by Wildberger \cite{{\bf 22}}.)

\head \S\sectg.  Loop algebras and periodic heaps \endhead

Having concentrated on the case of Kac--Moody algebras of finite type, we 
now turn our attention to the corresponding affine algebras.  For this 
purpose, it is convenient to introduce the notion of a periodic heap.

\definition{Definition \sectg.1}
Let $(E, \leq, \e)$ be a locally finite labelled heap over a graph $\Gamma$.
We call the labelled heap $(E, \leq, \e)$, and the associated heap
$[E, \leq, \e]$ {\it periodic} if there
exists a nonidentity automorphism $\phi : E \ra E$ of labelled posets
such that $\phi(x) \geq x$ for all $x \in E$.
\enddefinition

\remark{Remark \sectg.2}
\item{(i)}{It is immediate that any periodic heap is necessarily infinite.}
\item{(ii)}{The automorphism $\phi$ above restricts to an automorphism of 
the chains
$\e^{-1}(p)$ for $p$ a vertex of $\Gamma$.  By Remark \secta.7, this
automorphism must be of the form $\phi(E(p, x)) = E(p, x + t_p)$ for some
nonnegative integer $t_p$ depending on $p$ but not on the labelling chosen
for $E$, and furthermore, the 
automorphism $\phi$ can be reconstructed from the integers $t_p$.  If
$\a \in R^+$ is such that $\a(p) = t_p$, we will say that $\phi$ is
{\it periodic with period $\a$}.  If there is no automorphism $\phi'$ of $E$
with period $\a'$ such that $\a = n \a'$ with $n > 1$, then we also say that
$E$ is periodic with period $\a$.}
\endremark

\example{Example \sectg.3}
In the notation of Example \sectb.4, the heap $E$ is periodic with period 
$\cha(F)$.
\endexample

\proclaim{Lemma \sectg.4}
Let $E$ be a full heap over the Dynkin diagram $\Gamma$ of $A$, 
where $A$ is a simply laced generalized Cartan matrix of untwisted 
affine type (see the Appendix for examples, or \cite{{\bf 4}, Appendix}, 
\cite{{\bf 12}, \S4.8} for a 
complete list).  Let $A_0$ be the generalized Cartan matrix of finite type
obtained by omitting the row and column of $A$ corresponding to
the root $\a_0$.  Let $\d$ be the smallest positive imaginary root associated
to $A$.
\item{\rm (i)}{The root $\d$ is representable in $E$.}
\item{\rm (ii)}{The heap $E$ is periodic with period $\d$.}
\item{\rm (iii)}{Every positive root is representable in $E$.}
\endproclaim

\demo{Proof}
From \cite{{\bf 12}, Theorem 5.6 (b)}, we see that $\d = \th + \a_0$, where
$\th$ is the highest root of $\a$.  By \cite{{\bf 12}, Theorem 4.8 (c)}, 
$\lrt{\d}{\a_p^\vee} = 0$ for all vertices $p$, and since 
$\lrt{\a_0}{\a_0^\vee} = 2$, we have
$\lrt{\th}{\a_0^\vee} = -2$.  Let $L_\th \in \L_\th$; this exists by 
Proposition \secte.4 (i).  By Proposition \secte.4 (iii), we have
$[H_0, X_\th] = -2 X_\th$, which means that there exist vertices $\be, \be'$
of $E \backslash L_\th$, both labelled $0$, such that $L^+ = L_\th \cup \{\be\}
\cup \{\be'\}$ lies in $\L_{\d + \a_0}$ and such that $\be$ (respectively,
$\be'$) is a maximal (respectively, minimal) element of $L^+$.  By removing
either $\be$ or $\be'$ from $L$, we obtain an element of $\L_\d$, thus
proving (i).

For (ii), let $L_0 \in \L_\d$, which is a nonempty set by (i).  By \cite{{\bf 12},
Theorem 4.8 (c)}, we see that $L_0$ contains at least one element with each
possible label from $\Gamma$.  Let $\a$ be a maximal element of $L_0$, and 
let $p = \e(\a)$.  Since $\lrt{\d}{\a_p^\vee} = 0$, we see from Lemma \secte.3
(iii) that we can convert $L_0$ into an element $L'$ of $\L_\d$ 
(with $L' \ne L_0$)
by removing $\a$ and replacing it by a new minimal vertex with label $i$.
By repeating this procedure once for each element of the original heap $L_0$,
we obtain a heap $L_1 \in \L_\d$ such that $L_0 \cong L_1$ as heaps, 
$L_0 \cap L_1 = \emptyset$ and $L_0 \cup L_1 = L_1 \circ L_0$ is convex.  

By applying the above construction again to $L_1$, we obtain a sequence
$\{L_i \}_{i \geq 0}$ of disjoint isomorphic heaps $L_i \in \L_\d$ such that
$L_i \cup L_{i+1}$ is convex.  Since $L$ contains at least one element with 
each label, any element $\a$ with $\a \leq \be$ for some $\be \in L$
lies in one of the heaps $L_i$ for $i \geq 0$.

By a dual argument, we can also find a sequence $\{L_i\}_{i \leq 0}$ with
analogous properties, such that any element $\a$ with $\a \geq \be$ for some
$\be \in L$ lies in some $L_i$ for $i \leq 0$.  
Since $L_0$ contains an element with each possible label, the heaps $L_k$
for $k \in \zed$ partition the set $E$.  It follows that there is a
nonidentity automorphism $\phi : E \ra E$ of labelled posets sending $L_k$
to $L_{k-1}$ for any $k \in \zed$.  This construction shows that $\phi$ has
period $\d$.  Since $\d \ne n \d'$ for any $n > 1$ with $\d'
\in R^+$ (see \cite{{\bf 12}, Theorem 4.8 (c)}), we see that $E$ also has period
$\d$, completing the proof of (ii).

We next prove that if $\g$ is a root of $A_0$, then
$\d - \g$ is representable.  By Proposition \secte.4 (i), there exists
$L \in \L_\g$.  By Lemma \sectb.1 (vii), there is a proper ideal $I$
with $L \prec I$, and by periodicity, there exists $L' \in \L_\d$ with
$\phi^{-1}(I) = I \backslash L'$.  In this case, $L$ is a filter of
$L'$ and $L' \backslash L \in \L_{\d - \g}$, as required.

Part (iii) follows from this and
the fact (see \cite{{\bf 12}, \S7.4}) that the roots of $A$
are precisely the elements of $R^+$ of the form $\g + j \d$, where 
$j \in \zed$ and $\g$ is a root of $A_0$.  (The root is positive if 
either $\g > 0$ and $j \geq 0$, or $\g < 0$ and $j \geq 1$.)
\qed\enddemo

\proclaim{Lemma \sectg.5}
Let $E, \Gamma, A, \mu, f, \overline{A}, \overline{\Gamma}$ and $\overline{E}$
be as in Proposition \sectf.1, and assume in addition that $\overline{A}$ is
a generalized Cartan matrix of untwisted affine type.
Let $\d$ and $\bar{\d}$ be the smallest 
positive imaginary roots associated to $A$ and $\overline{A}$, respectively,
and suppose that we have $f(\d) = \bar{\d}$.
\item{\rm (i)}{The heap $\overline{E}$ is periodic with period $\bar{\d}$.}
\item{\rm (ii)}{Every positive root is representable in $\overline{E}$.}
\endproclaim

\demo{Proof}
Let $A_0$ (respectively, $\overline{A}_0$) be the generalized Cartan matrix
of finite type obtained by removing the zeroth row and column of $A$ 
(respectively, $\overline{A}$).

Since $\overline{A}$ is of untwisted affine type, it follows by \cite{{\bf 12},
\S7.4} that each positive root of $\overline{A}$ is of the form 
$n\bar{\d} + \g$, 
where $\g$ is a root for $\overline{A}_0$.  By construction of $\overline{A}$,
we have $\g = f(\g')$ for some root $\g'$ of $A_0$, and since we have
$f(\d) = \bar{\d}$ by hypothesis, it follows that $f(n\d + \g') = n\bar{\d}
+ \g$; note that $n\d + \g'$ is positive because $f$ respects positive and
negative roots.  Since $n\d + \g'$ is representable in $E$ by 
Lemma \sectg.4 (iii), part (ii) follows by folding.

Since $E$ is periodic with period $\d$, it follows that $\overline{E}$ is
periodic with period $\d'$, where $f(\d) = \bar{\d} = k \d'$ for some positive 
integer $k$.  Writing $\bar{\d} = \sum a_i \bar{\a_i}$, we have 
by \cite{{\bf 12}, Theorem 4.8 (c)} that
the $a_i$ are relatively prime integers, so we must have $k = 1$ and 
$\d' = \bar{\d}$, as required.
\qed\enddemo

\definition{Definition \sectg.6}
The automorphism $\phi$ of Lemma \sectg.4 (ii) and 
Lemma \sectg.5 (i) induces a permutation
(also denoted $\phi$) of the proper ideals of $E$.  We define an invertible
linear map $T : V_E \ra V_E$ by $T = X_\d$; note that $T^{-1} = Y_\d$.

If $\a_0$ is the simple root of $A$ such that $\d = \th + \a_0$,
then we define the {\it height}, $h(I)$ of a proper ideal $I$ of $E$ to be the 
maximum integer $t$ such that $E(0, t) \in I$.

We define the linear map $D : V_E \ra V_E$ by $D(v_I) = h(I) v_I$.
\enddefinition

\definition{Definition \sectg.7 (\cite{{\bf 12}, \S7})}
Let $\fg_0$ be a Kac--Moody algebra of finite type over $k = \complex$.
Then the {\it loop algebra} $\L(\fg_0)$ of $\fg_0$ is defined to be the $k$-vector
space $k[t, t^{-1}] \otimes_k \fg_0$ with Lie bracket defined by $
[P \otimes x, Q \otimes y] := PQ \otimes [x, y]
.$
\enddefinition

As explained in \cite{{\bf 12}, \S7},
a fundamental property of the untwisted affine Kac--Moody algebras
over $k = \complex$ is that
they can be constructed from loop algebras by adding both a one-dimensional 
centre $\complex K$ (using a universal central extension) and an additional
derivation $d$.  More precisely, we have

\proclaim{Theorem \sectg.8 (see \cite{{\bf 12}, \S7})}
Let $A$ be a generalized Cartan matrix of untwisted affine type
Let $A_0$ be the corresponding generalized Cartan matrix of finite type, with
Lie algebra $\fg_0$ over $k = \complex$ and highest root $\th$.  Then the 
Kac--Moody algebra $\fg$ associated to $A$ is the 
vector space $$
\L(\fg_0) \oplus \complex K \oplus \complex d
$$ equipped with the Lie bracket $$\eqalign{
[(t^m \otimes x) + \l K + \mu d, (t^n \otimes y) + \l_1 K + \mu_1 d]
=&
(t^{m+n} \otimes [x, y]) + \mu n (t^n \otimes y) \cr
&- \mu_1 m (t^m \otimes x) + m \d_{m, -n} (x | y) K, \cr
}$$ for $x, y \in \fg_0, \l, \mu, \l_1, \mu_1 \in \complex, m, n \in \zed$.
(See \cite{{\bf 12}, \S2.2} for the definition of $(x | y) \in \complex$.)
Given $\epsilon = \pm 1$, the isomorphism may be chosen to identify the 
subalgebra $\fg_0$ with
the subset $1 \otimes \fg_0$ of $\L(\fg_0)$, and to send the Chevalley generators
$e_0$ and $f_0$ of $\fg$ to $\epsilon t \otimes E_\th$ and 
$\epsilon t^{-1} \otimes -E_{-\th}$ respectively, where the $E_\a$ are as 
given in \cite{{\bf 12}, (7.8.5)}.
\endproclaim

\demo{Proof}
The complete argument may be found in \cite{{\bf 12}, \S7}.  The assertion about
the $\pm E_{\pm \th}$ follows from the fact (\cite{{\bf 12}, Remark 7.9 (c)})
that $\epsilon E_{\th}$ and $\epsilon E_{-\th}$ are exchanged by the Chevalley involution
of $\fg_0$, and that $$
[- \epsilon E_{-\th}, \epsilon E_{\th}] = - \th^\vee
$$ (see \cite{{\bf 12}, (7.8.5), (7.9.3)}).
\qed\enddemo

\proclaim{Lemma \sectg.9}
Let $A$ be a generalized Cartan matrix of untwisted affine type
with Dynkin diagram $\Gamma$, and let $E$ be a simply folded full heap over
$\Gamma$.
Let $A_0$ be the generalized Cartan matrix of finite type
obtained by omitting the row and column of $A$ corresponding to
the root $\a_0$, and let $\fg_0$ be the Lie algebra of $A_0$, identified with
the Lie algebra of operators $\fg_E$ on $V_E$ by Theorem \sectf.7.
\item{\rm (i)}{The map $T : V_E \ra V_E$ commutes with the operators
$X_L$, $Y_L$ and $H_L$.}
\item{\rm (ii)}{The loop algebra $\L(\fg_E)$ over $\complex$ acts faithfully 
on $V_E$ via $
(t^j \otimes P)(v_I) = T^j \circ P(v_I).
$}
\item{\rm (iii)}{We have $[D, T^j \circ P] = j T^j \circ P$ in the Lie algebra
of operators $\L(\fg_E) \oplus k D$.}
\endproclaim

\demo{Proof}
For part (i), we see that $T$ and $T^{-1}$ commute with $H_L$ by the fact
that $E$ is periodic of period $\d$ (Lemma \sectg.4 (ii) and
Lemma \sectg.5 (i)).
Again, by periodicity, we have $X_L . v_I \ne 0$ if
and only if $T \circ X_L . v_I \ne 0$ if and only if $X_L \circ T . v_I
\ne 0$.  By Lemma \sectd.6, we see that $$
T \circ X_L = \sgn(\d, \th) X_{\d + \a}
$$ and $$
X_L \circ T = \sgn(\th, \d) X_{\d + \a}
.$$  By \cite{{\bf 12}, (6.2.4), (7.8.3)}, we have $\sgn(\d, \th) = \sgn(\th, \d)$,
and it follows that $T$ (and therefore $T^{-1}$) commutes with the $X_L$.
A similar argument shows that $$
T^{-1} \circ Y_L = Y_L \circ T^{-1} = \sgn(\d, \th) Y_{\d + \a}
,$$ establishing the claim for the $Y_L$.

To prove (ii), we note that $P(v_I)$ is a linear combination of basis elements
$v_J$ with $h(J) = h(I)$ (where $h$ is as in Definition \sectg.6), but
that $T^j(v_I) = \pm v_{\phi^j(I)}$, where $h(\phi^j(I)) = h(I) + j$,
because $\a_0$ occurs in $\d$ with coefficient 1.  (The latter follows from
\cite{{\bf 4}, Proposition 17.2 (ii)}.)
Since $\fg_0$ is simple
and its action on $V_E$ is nontrivial, $\fg_0$ acts faithfully on $V_E$.
Part (ii) now follows from the fact that $T$ has infinite order.

Part (iii) follows from the observations that 
$D \circ T^j(v_I) = (h(I) + j)T^j(v_I)$ and
$T^j \circ D(v_I) = h(I) T^j(v_I)$.
\qed\enddemo

We can now state the main result of this section.

\proclaim{Theorem \sectg.10}
Let $A$ be a 
generalized Cartan matrix of untwisted affine type with Dynkin diagram
$\Gamma$, and let $\fg$ be the corresponding Kac--Moody algebra.  
Let $E$ be a simply folded full heap over $\Gamma$.
Let $\frak v$ be the Lie algebra of linear operators on $V_E$ over
$k = \complex$.
Then there is a homomorphism of Lie algebras $\psi : \fg \ra
{\frak v}$ sending the Chevalley generators $e_i, f_i$ (for $0 \leq i \leq n$)
to the generators $X_i, Y_i$ respectively, and sending the derivation $d$
to the operator $D$.  The kernel of $\psi$ is precisely $\complex K$.
\endproclaim

\demo{Proof}
Comparing the explicit formula of Theorem \sectg.8 (with $K$ acting as
zero) with the Lie algebra $\L(\fg_E) \oplus \complex D$ of Lemma \sectg.9,
and identifying $\fg_0$ with $\fg_E$ as in Theorem \sectf.7, 
we obtain a homomorphism from $\fg$ to 
$\L(\fg_E) \oplus \complex D$ with the required kernel, sending the
generators $e_i, f_i$ (for $0 < i \leq n$) to $X_i$ and $Y_i$ respectively.

To complete the proof, it is enough to show that we have $$\eqalign{
X_0 &= \psi(\epsilon T^{-1} \circ -E_{-\th}),\cr
Y_0 &= \psi(\epsilon T \circ E_{\th}),\cr
}$$ for some fixed $\epsilon = \pm 1$.  By Theorem \sectf.7 and Proposition 
\secte.4, we know that $\psi(E_\th) = X_\th$ and $\psi(-E_{-\th}) = Y_\th$.
By lemmas \sectg.4 and \sectg.5, $E$ is periodic with period $\d$, and the
fact that $\d = \th + \a_0$ then implies that we have $Y_\th . v_I \ne 0$
if and only if $X_0 . v_I \ne 0$ if and only if $T \circ Y_\th . v_I \ne 0$.
Now suppose that $Y_\th . v_I \ne 0$, define $L \in \L_\d$ to be such that
$\phi^{-1}(I) \cup L = I$, and define $L_\th \in \L_\th$ to be such that
$I = I' \cup L_\th$ with $Y_\th . v_I = v_{I'}$.  This means that 
$L = L_\th \cup \{\a\}$, where $\a$ is a minimal element of $L$ with label $0$.
By Lemma \sectd.4, we find that $$
\sgn(L) = \sgn(\{\a\}) \sgn(L_\th) \sgn(\th, \a) = \sgn(L_\th) \sgn(\th, \a)
.$$  The assertion for $X_0$ then follows by defining $\epsilon = 
\sgn(\th, \a)$.  The assertion for $Y_0$ is proved similarly, using 
$\epsilon' = \sgn(\a, \th)$.  By \cite{{\bf 12}, (6.2.4), (7.8.1), (7.8.3)}, 
it follows that $\epsilon = \epsilon'$, as required.
\qed\enddemo

\head \S\secth.  Quantum affine algebras, crystals, and the Weyl group 
action \endhead

Let $A$ be an $(l+1) \times (l+1)$ generalized Cartan matrix of affine type
corresponding to an untwisted Kac--Moody algebra $\fg$ over $k = \complex$.
We assume that $A$ is indexed in such a way that the removal of the zeroth row
and column of $A$ results in the generalized Cartan matrix for the 
corresponding Kac--Moody algebra of finite type.  The sets $\Pi$ and $\Pi^\vee$
were defined in \S\sectc\  for the set $I = \{0, 1, \ldots, l\}$.  
We now extend $I$ by redefining it as $\{0, 1, \ldots, l+1\}$, and we
extend the sets 
$\Pi$ and $\Pi^\vee$ accordingly.  The $\zed$-bilinear form is then extended
by setting $$\eqalign{
\lrt{\a_0}{\a_{l+1}^\vee} &= 1,\cr
\lrt{\a_i}{\a_{l+1}^\vee} &= 0 \text{ if } 1 \leq i \leq l,\cr
\lrt{\a_{l+1}}{\a_0^\vee} &= 1,\cr
\lrt{\a_{l+1}}{\a_i^\vee} &= 0 \text{ if } 1 \leq i \leq l,\cr
\lrt{\a_{l+1}}{\a_{l+1}^\vee} &= 0;\cr
}$$ this corresponds to the action of $H^*$ on $H$ described in 
\cite{{\bf 4}, \S17.1}.

There are several slight variants in the literature of the definition of 
a quantized affine algebra; ours is based on the one in 
\cite{{\bf 2}}.  Let $q$ be an indeterminate. For
nonnegative integers $n\ge r$, define $$\eqalign{
[n] &= {{q^n - q^{-n}} \over {q - q^{-1}}},\cr
[n]! &= \prod_{i = 1}^n [i],\cr
\qchoose{n}{r} &= {{[n]!} \over {[r]! [n-r]!}}.\cr
}$$  Since we are in the untwisted case, \cite{{\bf 4}, Proposition 17.2 (ii)}
shows that, in the notation of \cite{{\bf 2}, \S2.2}, we have $d = 1$.  By
\cite{{\bf 4}, Proposition 17.9 (a, b)}, the symbol $q_i$ in \cite{{\bf 2}, \S2.2}
is equal to $q$ if $\a_i$ is a short root, and to $q^2$ if $\a_i$ is a long
root.

\definition{Definition \secth.1}
Define the quantum affine algebra $\U$ to be the associative algebra
with $1$ over $\kyu(q)$ generated by elements $E_i$, $F_i$
($i\in I$), $q^h$ (for $h\in \zed\Pi^\vee$),
with defining relations $$\eqalign{
q^0 &= 1,\cr
q^h q^{h'} &= q^{h+h'},\cr
q^h E_i q^{-h} &= q^{\lrt{\a_i}{h}} E_i,\cr
q^h F_i q^{-h} &= q^{-\lrt{\a_i}{h}} F_i,\cr
E_i F_j - F_j E_i &= \delta_{ij}
{{t_i - t_i^{-1}} \over {q_i - q_i^{-1}}}, \cr
\sum_{p=0}^{b}(-1)^p E_i^{(p)} E_j E_i^{(b-p)} &= 
\sum_{p=0}^{b}(-1)^p F_i^{(p)} F_j F_i^{(b-p)} = 0 \quad
\text{for $i\ne j$,}
}$$ where $q_i$ is as above,
$t_i = q_i^h$,
$b = 1 - \lrt{\a_j}{\a_i^\vee}$,
$E_i^{(p)} = E_i^p/[p]!$, and
$F_i^{(p)} = F_i^p/[p]!$.
\enddefinition

The {\it weight lattice} $P$ of $\U$ is the $\zed$-module 
$\Hom_\zed(\zed\Pi^\vee, \zed)$.  Since $P$ and $\zed\Pi^\vee$ are free
$\zed$-modules of the same finite rank, they are in natural duality.

A $\U$-module $M$ is called {\it integrable\/} if
\item{(a)}{all $E_i$, $F_i$ ($i\in I$) act locally nilpotently; that
is, for each $v \in M$ we have $E_i^N . v = F_i^N . v = 0$ for sufficiently
large $N$, and}
\item{(b)}{$M$ admits a {\it weight space decomposition}: $$
   M = \bigoplus_{\lambda \in P} M_\lambda,
   \text{\ where\ } M_\lambda = \{ u\in M : q^h u
     = q^{\lrt{\lambda}{h}} u \text{\ for all\ }
h \in \zed\Pi^\vee\}.$$}

Let $A_0$ be the subring of $\kyu(q)$ consisting of the
rational functions of $q$ that are regular at $q=0$.

Let $M$ be an integrable $\U$ module.  Kashiwara \cite{{\bf 13}} showed that
we have $$
M=\bigoplus_{\lambda} F_i^{(n)}(\ker E_i\cap M_\lambda)
,$$ and defined operators $\kase_i, \kasf_i : M \ra M$ for each
$0 \leq i \leq l$ (often called
{\it Kashiwara operators}) by $$
\kasf_i(F_i^{(n)}u)=F_i^{(n+1)}u\text{ and
}\kase_i(F_i^{(n)}u)=F_i^{(n-1)}u,$$ where $u \in \ker E_i\cap M_{\lambda}$.
(We interpret $F_i^{(-1)}u = 0$ above.)
Any element of $M$ is uniquely expressible as a sum of such elements 
$F_i^{(n)}u$.

\definition{Definition \secth.2}
Let $M$ be an integrable $\U$-module. A pair $(\L,\Be)$ is called a
{\it crystal basis} of $M$ if it satisfies:
\item{\rm (i)}{$\L$ is a free $A_0$-submodule of $M$ such that $M\cong
\kyu(q)\otimes_{A_0}\L$,}
\item{\rm (ii)}{$\L = \bigoplus_{\lambda\in P} \L_\lambda$ where
$\L_\lambda = \L\cap M_\lambda$ for $\lambda\in P$,}
\item{\rm (iii)}{$\Be$ is a $\kyu$-basis of $\L/q\L \cong \kyu 
\otimes_{A_0}\L$,}
\item{\rm (iv)}{$\kase_i\L\subset\L$, $\kasf_i\L\subset\L$ for all $i\in
I$,}
\item{\rm (v)}{if we denote operators on $\L/q\L$ induced by $\kase_i$, $
\kasf_i$ by the same symbols, we have
$\kase_i\Be\subset\Be\sqcup\{0\}$, $\kasf_i\Be\subset\Be\sqcup\{0\}$,}
\item{\rm (vi)}{for any $b,b'\in\Be$ and $i\in I$, we have $b' = \kasf_i b$ 
if and only if $b = \kase_i b'$.}
\enddefinition

The following result is a $q$-analogue of Theorem \sectg.10 (but without
the assertion about the kernel).

\proclaim{Theorem \secth.3}
Let $A$ be 
generalized $(l+1) \times (l+1)$ Cartan matrix of untwisted affine type
and let $\fg$ be the corresponding Kac--Moody algebra.  
Let $E$ be a simply folded full heap over $\Gamma$.
\item{\rm (i)}{Over the field $k = \kyu(q)$, 
$V_E$ has the structure of a (left) $\U$-module
such that $E_i$ and $F_i$ act as $X_i$ 
and $Y_i$, and such that for all $h \in \zed\Pi^\vee$ we have 
$q^h . v_I = q^\l v_I$, where $\l \in \zed$ is such that $h . v_I = \l v_I$.
Here, we identify $\a_i^\vee$ with $H_i$ for $0 \leq i \leq l$, and 
with the operator $D$ for $i = l+1$.}
\item{\rm (ii)}{The $\U$-module $V_E$ is integrable and has a crystal
basis $(\L, \Be)$, where $$
\Be = \{v_I : I \text{ is a proper ideal of } E \}
$$ and $\L$ is the free $A_0$-module with $\Be$ as a basis.}
\endproclaim

\demo{Proof}
Apart from the relations involving $q^{\a_{l+1}^\vee}$ 
and $D$, part (i) follows by imitating the
proof of Theorem \sectc.1, substituting exponentials where necessary. 
It follows from the definitions that $[D, X_i] = \d_{0i} X_i$, and that
$[D, Y_i] = -\d_{0i} Y_i$.  The remaining assertions of (i) now follow from the
definition of $\lrt{\a_i}{\a_{l+1}^\vee}$.

The action of $E_i$ and $F_i$ on $V_E$ is locally nilpotent by Lemma
\sectb.7 (7).  Since $H_i(v_I)$ and $D(v_I)$ are integer multiples 
of $v_I$, it follows that for $h \in \zed\Pi^\vee$, we have $$
q^h . v_I = q^{\l(h)} v_I
$$ for some $\l \in P$ depending on $I$ but not $h$.  This shows that each
$v_I$ is a weight vector.  Parts (i), (ii) and (iii) of Definition \secth.2
now follow from the definition of $\L$.

Now fix generators $E_i$ and $F_i$, and a basis element $v_I$.  Since 
$E_i^2$ acts as zero, and the action of $E_i$ and $F_i$ takes basis elements
either to other basis elements, or to zero, we see that either $v_I$ lies
in $\ker E_i$, or that $v_I = F_i . v_{I'}$ for some $v_{I'} \in \ker E_i$.
Using the definition of raising and lowering operators, we now find that 
the Kashiwara operators $\kasf_i$ and $\kase_i$ are simply given by the
actions of $F_i$ and $E_i$ respectively.  Parts (iv), (v) and (vi)
of Definition \secth.2 now follow.
\qed\enddemo

\remark{Remark \secth.4}
Notice that the distributive lattice structure induced by Corollary \sectb.2
is compatible with the partial order induced on the basis by the Kashiwara
operators.  It would be interesting to know whether this phenomenon is 
typical.

Much of the literature about crystals deals with the case of crystals
with extremal weight vectors, but the crystals mentioned in the theorem
do not have this property.  Another possible approach to these crystals would
be to bypass the quantum affine algebra and use Kashiwara's notion of {\it
abstract crystals}, which are crystals equipped with formal weight functions, 
satisfying certain axioms.  We do not pursue this here for reasons of space.
\endremark

\proclaim{Proposition \secth.5}
Maintain the assumptions of Theorem \secth.3, and assume that $\Gamma$
has finitely many vertices.  Then
the $\U$-module $V_E$ is cyclic and is generated by any one basis element
$v_I$.
\endproclaim

\demo{Proof}
Let $v_I$ and $v_J$ be basis elements.  It is enough to exhibit an element 
$u \in \U$ such that $u . v_I = v_J$.  Let $K = I \cap J$.  Since 
$K \subseteq I$, $I \backslash K$ is finite by Lemma \sectb.1 (vi), 
so it follows that there is a finite sequence $j_1, \ldots, j_l$ such that $$
F_{j_1} F_{j_2} \cdots F_{j_l} . v_I = v_K
.$$  A similar argument shows that there is a finite sequence $i_1, \ldots,
i_k$ such that $$
E_{i_1} E_{i_2} \cdots E_{i_k} . v_K = v_J
.$$  Concatenating these sequences produces the required element $u$.
\qed\enddemo

In \cite{{\bf 14}}, Kashiwara introduces the notion of a ``normal crystal'' (now
often referred to as a ``regular crystal''; see \cite{{\bf 2}, \S2.8}).  Such 
crystals naturally
carry an action of the associated Weyl group.  In this section, we show
that the crystals of Theorem \secth.3 also have this property, and that
furthermore, the action factors through a Temperley--Lieb type quotient.
This brings to light some representation theoretic obstructions to finding
full heaps for certain Kac--Moody algebras.

\definition{Definition \secth.6}
Let $A$ be a generalized Cartan matrix 
with Dynkin diagram $\Gamma$, and let $E$ be a simply folded full heap over
$\Gamma$.
For each $i \in I$, we define a linear operator $S_i$ on $V_E$
by requiring that $$
S_i(v_I) = \cases
F_i(v_I) & \quad \text{if } F_i(v_I) \ne 0,\cr
E_i(v_I) & \quad \text{if } E_i(v_I) \ne 0,\cr
v_I & \quad \text{otherwise.}
\endcases$$
\enddefinition

The definition of full heap guarantees that the cases in the above definition
do not overlap.

\proclaim{Proposition \secth.7}
Suppose that $A$, $E$ and $\Gamma$ satisfy the hypotheses of 
Definition \secth.6, and let $W = W(\Gamma)$ be the Weyl group of $\Gamma$.
\item{\rm (i)}{The assignment $s_i \mapsto -S_i$ defines a unique (left)
$kW$-module structure on $V_E$.}
\item{\rm (ii)}{If $s_i, s_j$ are a pair of noncommuting generators of $W$
and the subgroup of $W$ they generate, $W_{ij} = \lan s_i, s_j \ran$, is
finite, then the element $$
\sum_{w \in W_i} w
$$ of $kW$ annihilates $V_E$.}
\item{\rm (iii)}{The $kW$-module $V_E$ is cyclic, and any of the basis elements
$v_I$ is a generator.}
\endproclaim

\demo{Proof}
To prove (i), we need to check the defining relations of the Weyl group.
The relation $s_i^2 = 1$ holds by Theorem \secth.3 (ii) and Definition
\secth.2 (vi).  If $a_{ij} = 0$, the relation $s_i s_j = s_j s_i$ follows 
from the fact that each element of $\{E_i, F_i\}$ commutes with each element
of $\{E_j, F_j\}$.

Now let $s_i$ and $s_j$ be a pair of noncommuting generators, and let
$v_I$ be a basis element.  

Let us first suppose that either $F_i . v_I \ne 0$
or $F_j . v_I \ne 0$.  Since a heap cannot have two maximal vertices with
adjacent labels, these possibilities are mutually exclusive, so without
loss of generality, we may assume that $F_i . v_I \ne 0$.  We may now invoke
Lemma \sectb.5 with $p = i$ and $q = j$.

There are five subcases to consider.  In case 1 (respectively, 2), we have 
$s_i s_j s_i = s_j s_i s_j$ and case (i) (respectively, (ii)) of Lemma
\sectb.5 applies.  In case 3, (respectively, 4, 5) we have 
$s_i s_j s_i s_j = s_j s_i s_j s_i$ and case (i) (respectively, (ii), (iii))
of Lemma \sectb.5 applies.  In each case, the Weyl group relation is
respected by the claimed module action, and we have $
\sum_{w \in W_i} w . v_I = 0
.$  For example, in case 4, the identity and $S_i S_j S_i$ both act as the
identity; $S_i$ and $S_j S_i$ both act as $F_i$; $S_j$, $S_i S_j S_i S_j$
and $S_j S_i S_j S_i$ all act as $E_j$; and $S_i S_j$ and $S_j S_i S_j$ both
act as $E_i E_j$.

Now let us suppose that $F_i . v_I = F_j . v_I = 0$.  If we also have
$E_i . v_I = E_j . v_I = 0$, both $S_i$ and $S_j$ acts as the identity
on $v_I$, and it is clear that the Weyl group relation holds and that
$\sum_{w \in W_i} w . v_I = 0$, so we may assume that this is not the case.  
As before, the conditions $E_i . v_I \ne 0$ and $E_j . v_I \ne 0$ are
mutually exclusive, so we may assume that $E_i . v_I \ne 0$ without loss
of generality.  There must therefore exist a minimal element $\a \in E
\backslash I$ with $\e(\a) = i$.

In either case, we apply
Lemma \sectb.5 to the ideal $I \cup \{\a\}$, with $p = i$ and $q = j$.
There are two subcases to consider, according as $s_i s_j s_i = s_j s_i s_j$
or $s_i s_j s_i s_j = s_j s_i s_j s_i$.  
In either case, part (ii) of the lemma must apply.  For example, in the
former case, the identity and $S_j$ both act as the identity; $S_i$ 
and $S_i S_j$ both act as $E_i$; and $S_j S_i$, $S_i S_j S_i$ and
$S_j S_i S_j$ all act as $E_j E_i$.  We conclude that the Weyl group
relation holds and that $\sum_{w \in W_i} w . v_I = 0$.

This completes the proofs of parts (i) and (ii).
Part (iii) follows by the same argument used to prove Proposition \secth.5.
\qed\enddemo

\proclaim{Corollary \secth.8}
There are no full heaps over any Dynkin diagrams of finite type, or over
those of types $F_4^{(1)}$, $E_8^{(1)}$, or $E_6^{(2)}$.
\endproclaim

\demo{Proof}
Let $J$ be the ideal of $\zed W$ generated by the elements 
$\sum_{w \in W_i} w$ for each subgroup $W_i$ of $W$ generated by a pair
of noncommuting generators.  The algebra $\zed W/ J$ is precisely the 
generalized Temperley--Lieb algebra $TL(W)$ of \cite{{\bf 6}}, with the parameter
$q$ specialized to 1.  (The sign twists in Proposition \secth.7 (i) were
inserted for compatibility of the ideal $J$ with \cite{{\bf 6}}.)

If the hypotheses of Proposition \secth.7 hold, then we may set $k = \kyu$,
and the definition of full heap implies that $V_E$ is an infinite 
dimensional cyclic $\kyu \otimes_\zed TL(W)$-module.  This means that
$TL(W)$ and its $q$-analogue must have infinite rank.  This cannot
happen if $A$ is of finite type, because in this case $W$ is well known to
be a finite group.  If $A$ is of type $F_4^{(1)}$ (respectively, $E_8^{(1)}$, 
$E_6^{(2)}$) then the algebra $TL(W)$ is of type $F_5$ (respectively, $E_9$,
$F_5$) in the notation of \cite{{\bf 6}}, and is of finite rank by \cite{{\bf 6},
Theorem 7.1}.
\qed\enddemo

\head \S\secti. Applications and questions \endhead

We now outline how the results of this paper can be used to simplify those 
described by Wildberger \cite{Wi} and generalize them to the non simply laced 
case.  

Let $E$ be a simply folded heap over the Dynkin diagram $\Gamma$ of an 
untwisted affine Kac--Moody algebra, where $\Gamma$ is labelled so 
that vertex $p_0$ is the additional vertex relative to the corresponding 
finite type algebra, $\fg_0$, and where $E$ is labelled as in Remark \secta.7.
Let $$
E' = \{ x \in E : x \leq E(p_0, 0) \}
,$$ and let $E_0$ be the subheap of $E$ consisting of the vertices $$
E_0 = \{ x \in E : x \not\geq E(p_0, 1) \} \backslash E'
,$$ where we regard $E'$ and $E_0$ as subheaps of $E$.  It is straightforward
to check that the map $J \mapsto J \cup E'$ defines an order-preserving
bijection between the ideals of $E_0$ and the proper ideals of $E$ of height
zero, so that the ideals of $E_0$ are in natural bijection with the orbits
of proper ideals of $E$ under the action of the automorphism $\phi$.
This leads to an irreducible representation of the simple Lie algebra
$\fg_0$ over $\complex$, where the highest and lowest weight vectors correspond
to the ideals $E_0$ and $\emptyset$ of $E_0$, respectively.

Wildberger's approach is to work directly with the set of ideals of $E_0$
(but only in the simply laced case), using raising and lowering operators 
similar to those in this paper.  Our approach is simpler in that (a) we do not
need to impose a partial order on the set of convex subheaps and 
(b) the parity of a root heap, in the simply laced case, can be easily 
computed intrinsically 
in terms of its isomorphism type, whereas in \cite{Wi} the parity of a 
root heap in $\L_\a$ is computed by comparing it with a canonical
representative of $\L_\a$.  The approach here gives us enough control over
the signs that we can make a precise link with the Chevalley bases in
\cite{{\bf 12}}, for each possible orientation of the diagram.

When we apply this technique to the simple Lie algebras of type $B$, we obtain
a combinatorial construction of the spin representation that does not involve
Clifford algebras, analogous to the constructions described by Wildberger
\cite{Wi, \S5} for type $D$.  (The reader is referred to \cite{{\bf 4}, \S13.5}
for the Clifford algebra construction.)  In these cases, the finite heap 
$E_0$, whose size grows quadratically with the rank, may be much 
smaller than the dimension of the spin module, which grows exponentially
with the rank.  The number of ideals of $E_0$ is a power of $2$ in this case,
which may be shown by exhibiting a bijection between ideals of $E_0$ and
certain paths; we hope to give details of this elsewhere.  It would be
interesting to see if the Clifford algebra itself has an action by raising
and lowering operators on the spaces $V_{E_0}$ and $V_E$.

It may be conjectured that the aforementioned representation of the finite 
dimensional Lie algebra $\fg_0$ will be a minuscule representation, and that
there will be a $1$--$1$ correspondence between simply folded full heaps 
over untwisted affine Dynkin diagrams and minuscule representations of
simple Lie algebras over $\complex$.  
Such a result would require a classification of full heaps over
untwisted affine Dynkin diagrams.  When the Dynkin diagram contains no
circuits, \idest in types other than $A$, this is relatively
easy because the heaps are ranked as posets; the latter may be proved
by \cite{{\bf 8}, Theorem 2.1.1 (iii)} and has an analogue for minuscule 
heaps \cite{{\bf 20}, Corollary 3.4}.  In type $A_l^{(1)}$, things are more 
complicated, but based on the results of \cite{{\bf 24}},
we expect that there will be $l$ isomorphism classes, most of which will 
not be ranked.  A good context to examine these might be
the extended slant lattices of Hagiwara \cite{{\bf 10}, \S8}.

The crystal bases for minuscule modules for simple Lie algebras have been
known for some time; for example, they are implicit (in the form of
canonical bases) in the work of
Lusztig \cite{{\bf 18}, Theorem 19.3.5, Proposition 28.1.4}.  Another
construction of these bases may be given by restricting the crystal basis 
$\Be$ arising from a simply folded full heap $E$ to the finite dimensional 
module corresponding to the finite subheap $E_0$.  An advantage of our
approach is that one can describe the action of a Chevalley basis 
on the canonical basis.

Full heaps also exist over certain infinite Dynkin diagrams, such as the
diagrams $A_\infty$, $B_\infty$ and $D_\infty$ of \cite{{\bf 12}, \S7.11}.  In
these cases, the heaps arising are reminiscent of those needed to construct
spin representations in the finite case.  However, we do not know any 
examples of full heaps over finite graphs that do not correspond to 
affine Kac--Moody algebras.

In a future paper, we will show that $V_E$ has an interesting structure
as a module for the affine Weyl group.  
In many cases, this gives new, uniform constructions of representations
familiar from other contexts: for example, in types $A_l^{(1)}$ 
and $E_6^{(1)}$ the module structure appears to agree with certain of the
author's cell modules for tabular algebras (see \cite{{\bf 7}, \S6} and 
\cite{{\bf 9}, \S1.2}) after specializing the parameter to $1$.

\head Acknowledgements \endhead

I am grateful to P.P. Martin for suggesting improvements to an early version
of this paper, and to R.J. Marsh for some very helpful correspondence.  I
also thank R.A. Liebler, J. Losonczy and G.E. Moorhouse for their encouraging
comments.

\head Appendix: Examples of simply folded full heaps \endhead

In the appendix, we give some examples of full heaps over Dynkin diagrams
of affine Kac--Moody algebras.  All these heaps are periodic, and the 
dashed boxes in the diagrams indicate the repeating motif.  (Note that, 
for an untwisted affine algebra, the number of elements in the dashed box
is the Coxeter number of the associated finite type algebra.)

Folding by each of the automorphisms $\mu$ shown below either leaves the 
period of a heap $E$ unchanged, or, in the case of a twisted affine algebra 
or type $A_1^{(1)}$, the period is halved.  In the latter case, the number 
of orbits of proper ideals of $E$ under the action of $\phi$ is also halved.

Recall that if $(E, \leq)$
is a partially ordered set, a function $\rho : E \ra \zed$ is said to be a 
{\it rank function} for $(E, \leq)$ if whenever $a, b \in E$ are such that
$a < b$ is a covering relation, we have $\rho(b) = \rho(a) + 1$.  If a
rank function for $(E, \leq)$ exists, we say $(E, \leq)$ is {\it ranked}.
The heaps shown in this section are all ranked.  

\vfill\eject

\subhead Type $A_l^{(1)} (l > 1)$, natural representation \endsubhead

The Dynkin diagram $\Gamma$ of type $A_l^{(1)}$ (for $l > 1$) 
is labelled as in Figure 8.

\topcaption{Figure 3} The Dynkin diagram of type $A_l^{(1)} (l > 1)$
\endcaption 
\centerline{
\hbox to 2.888in{
\vbox to 1.125in{\vfill
        \includegraphics{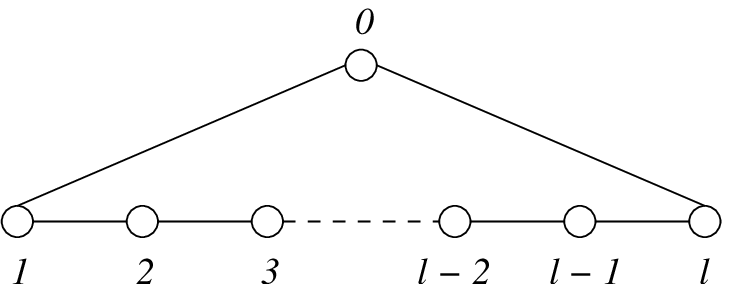}
}
\hfill}
}

The finite heap $E_0$ corresponding to the full heap $E$ shown in 
Figure 4 gives rise to the natural representation
of the simple Lie algebra of type $A_l$ (see \cite{{\bf 4}, \S8.1} for more 
details).  In Wildberger's notation \cite{{\bf 22}}, we have $E_0 = F(A_l, 1)$.
It is clear that the heap $E$ has $l+1$ orbits of proper ideals under $\phi$.

\topcaption{Figure 4} A full heap, $E$, over the Dynkin diagram of type 
$A_l^{(1)} (l > 1)$
\endcaption
\centerline{
\hbox to 2.402in{
\vbox to 2.402in{\vfill
        \includegraphics{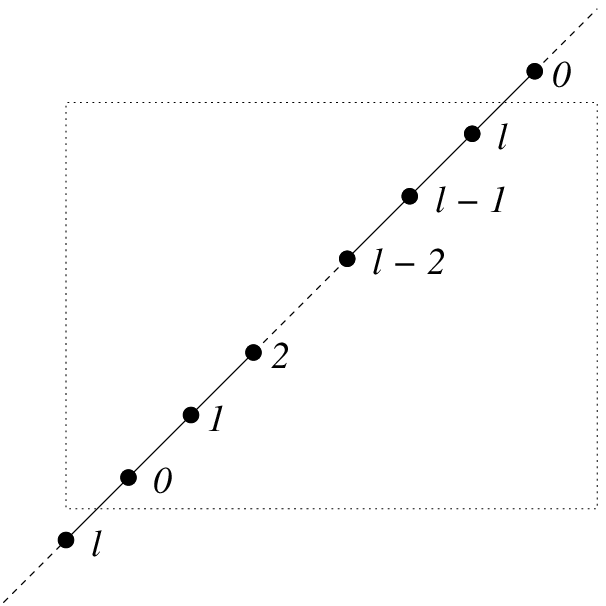}
}
\hfill}
}

The Dynkin diagram $\Gamma$ has an automorphism $\mu$ given by sending
vertex $i$ to vertex $l+1-i$ and fixing vertex $0$.  If $l$ 
is odd, then $\mu$ also fixes $(l+1)/2$, and we obtain the Dynkin diagram 
$\overline{\Gamma}$ of type $C_{l'}$, where $l = 2l' - 1$.

\topcaption{Figure 5} The Dynkin diagram of type $C_l^{(1)} (l > 1)$
\endcaption 
\centerline{
\hbox to 2.888in{
\vbox to 0.388in{\vfill
        \includegraphics{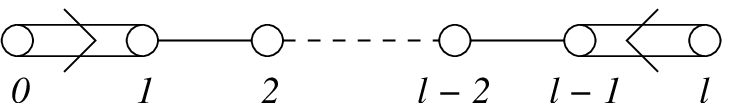}
}
\hfill}
}

The heap $E$ folds to a heap $\overline{E}$ over $\overline{\Gamma}$, via
$\mu$.  The heap $\overline{E}$ has $2l$ orbits of proper ideals.

\topcaption{Figure 6} A full heap, $\overline{E}$, over the Dynkin 
diagram of type $C_l^{(1)} (l > 1)$
\endcaption
\centerline{
\hbox to 2.152in{
\vbox to 2.777in{\vfill
        \includegraphics{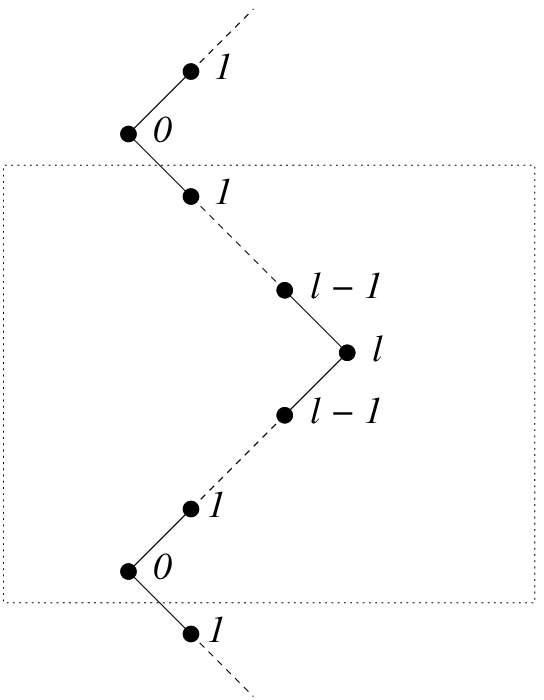}
}
\hfill}
}

The Dynkin diagram of type $A_1^{(1)}$ is as shown in Figure 7.

\topcaption{Figure 7} The Dynkin diagram of type $A_1^{(1)}$
\endcaption 
\centerline{
\hbox to 0.888in{
\vbox to 0.388in{\vfill
        \includegraphics{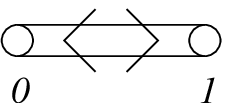}
}
\hfill}
}

If $l = 3$, so that the Dynkin diagram $\Gamma$ of type $A_l^{(1)}$ is a 
square, 
there is an automorphism $\mu$ of $\Gamma$ given by rotation by a half turn.
In this case, the full heap of Figure 4 folds to a full heap over the
Dynkin diagram of type $A_1^{(1)}$.  (This causes the period to halve, which
is behaviour usually characteristic of the twisted affine case.)  The results
of this paper may be checked by hand for this case.

\vfill\eject

\subhead Type $D_l^{(1)}$, natural representation \endsubhead

The Dynkin diagram $\Gamma$ of type $D_l^{(1)}$ is labelled as in Figure 8.

\topcaption{Figure 8} The Dynkin diagram of type $D_l^{(1)}$
\endcaption
\centerline{
\hbox to 3.486in{
\vbox to 1.138in{\vfill
        \includegraphics{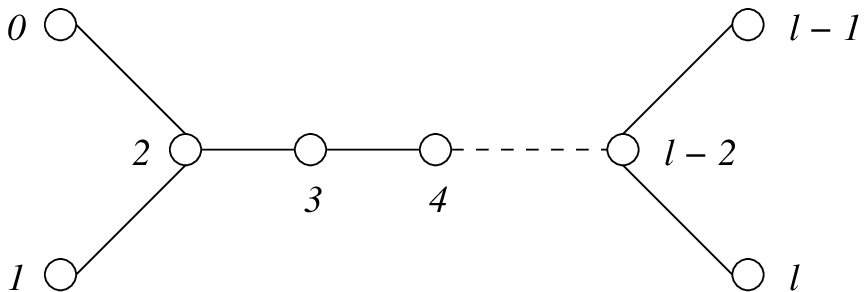}
}
\hfill}
}

The finite heap $E_0$ corresponding to the full heap $E$ shown in Figure 9 
gives rise to the natural representation
of the simple Lie algebra of type $D_l$; see \cite{{\bf 4}, \S8.2} for details
of another construction of this representation.  In Wildberger's notation,
we have $E_0 = F(D_l, l-1)$.  The heap $E$ has $2l$ orbits of proper ideals
under $\phi$.

\topcaption{Figure 9} A full heap, $E$, over the Dynkin diagram of type 
$D_l^{(1)}$
\endcaption
\centerline{
\hbox to 2.152in{
\vbox to 3.527in{\vfill
        \includegraphics{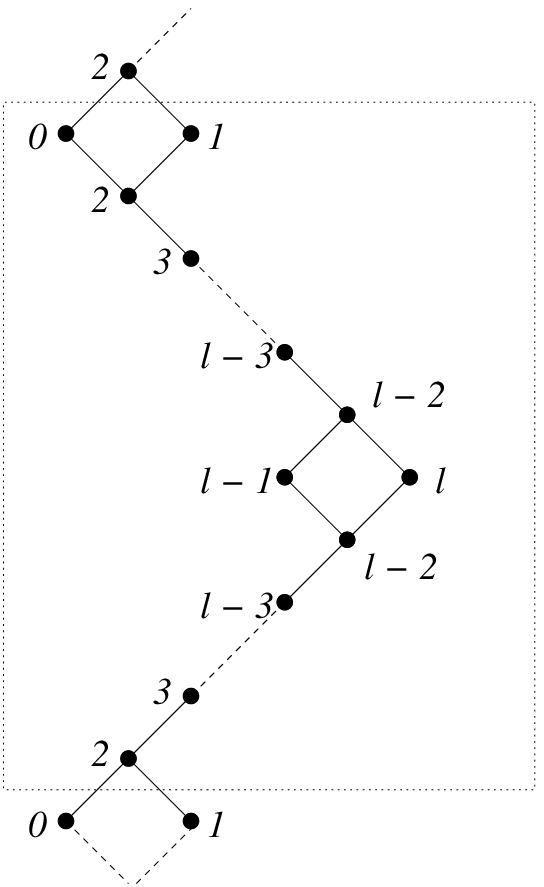}
}
\hfill}
}

The Dynkin diagram $\Gamma$ has an automorphism $\mu$ given by sending
vertex $i$ to vertex $l+1-i$.  If $l$ is odd, then $\mu$ has a fixed point
and we obtain the Dynkin diagram $\overline{\Gamma}$ of type 
$A_{2l-1}^{(2)}$ from
that of type $D_{2l-1}^{(1)}$; the former is shown in Figure 10.
The orbit $\{i, 2l-i\}$ in type $D_{2l-1}^{(1)}$ corresponds to the vertex $i$ 
in type $A_{2l-1}^{(2)}$.

\topcaption{Figure 10} The Dynkin diagram of type $A_{2l-1}^{(2)}$
\endcaption
\centerline{
\hbox to 3.569in{
\vbox to 1.138in{\vfill
        \includegraphics{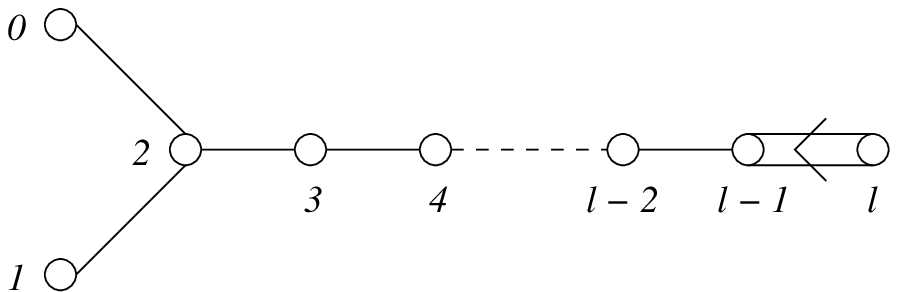}
}
\hfill}
}

The heap $E$ folds to a heap $\overline{E}$ over $\overline{\Gamma}$, via
$\mu$, shown in Figure 11.

\topcaption{Figure 11} A full heap, $\overline{E}$, over the Dynkin 
diagram of type $A_{2l-1}^{(2)}$
\endcaption
\centerline{
\hbox to 2.152in{
\vbox to 3.527in{\vfill
        \includegraphics{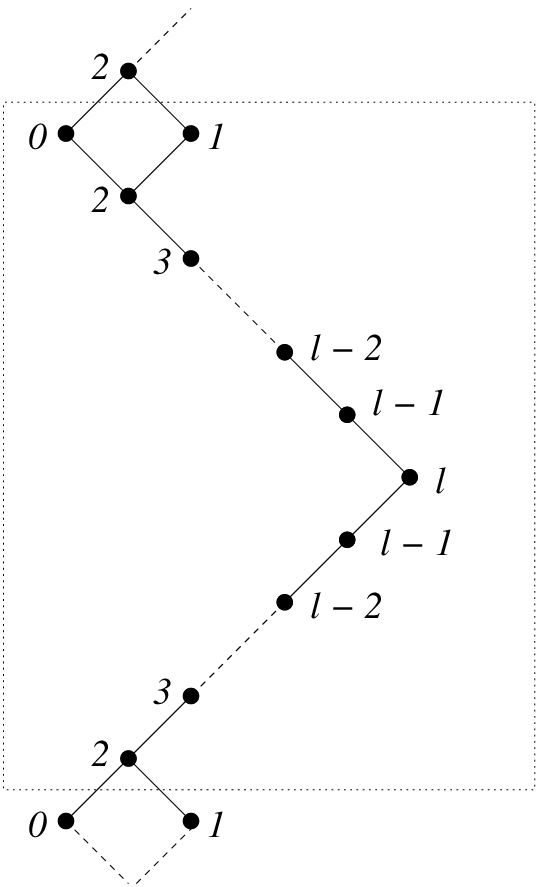}
}
\hfill}
}

\vfill\eject

\subhead Type $D_l^{(1)}$, spin representations \endsubhead

We now consider full heaps $E$ over the Dynkin diagram of type 
$D_l^{(1)}$ (see Figure 8) corresponding to a spin representation of the 
simple Lie algebra of type $D_l$.  The heap $E$ is ranked; we will call the 
subheap of $E$ given by the elements of rank $k$ the ``$k$-th layer'' of $E$.  
The even numbered vertices in the set $X = \{2, 3, 4, \ldots, l-2\}$ occur
in the $k$-th layer if and only if $k$ is even, and the odd numbered vertices
in $X$ occur in the $k$-th layer if and only if $k$ is odd.  
If $l$ is odd, then the $k$-th layer contains precisely one other vertex,
labelled $l-1$, $0$, $l$ or $1$ according as $k = 0$, $1$, $2$ or $3 \mod 4$.
If $l$ is even, then the $(2k+1)$-st layer contains precisely two other
vertices, whose labels are $0$ and $l-1$ if $k$ is odd, and $1$ and $l$ if
$k$ is even.  Figure 12 shows examples of such heaps for $l = 6$ and $l = 7$.
Another isomorphism class of heaps may be obtained in each case by twisting
by the graph automorphism exchanging vertices $l-1$ and $l$.  The corresponding
finite heaps $E_0$ are $F(D_l, 0)$ and $F(D_l, 1)$ in Wildberger's notation.
There are $2^{l-1}$ orbits of proper ideals of $E$ under $\phi$.

\topcaption{Figure 12} Full heaps over the Dynkin diagram of type 
$D_l^{(1)}$ for $l = 6$ and $l = 7$
\endcaption
\centerline{
\hbox to 4.152in{
\vbox to 1.777in{\vfill
        \includegraphics{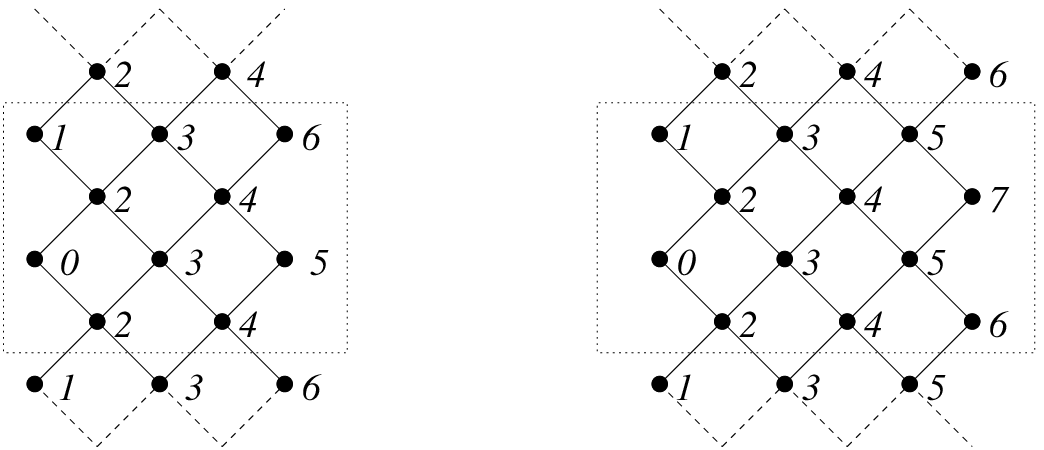}
}
\hfill}
}

There is an automorphism $\mu_1$ of the Dynkin diagram of $D_l^{(1)}$ obtained
by exchanging vertices $l-1$ and $l$, and fixing each of the other vertices.
The full heaps over $D_l^{(1)}$ discussed above fold via $\mu_1$ to a heap 
over the Dynkin
diagram of type $B_{l-1}^{(1)}$ (see Figure 13).  The only difference this 
makes to the full heap is that the vertices formerly numbered $l$ all
have their labels changed to $l-1$; this has the effect of merging the two
isomorphism classes.

\topcaption{Figure 13} The Dynkin diagram of type $B_l^{(1)}$
\endcaption
\centerline{
\hbox to 3.069in{
\vbox to 1.138in{\vfill
        \includegraphics{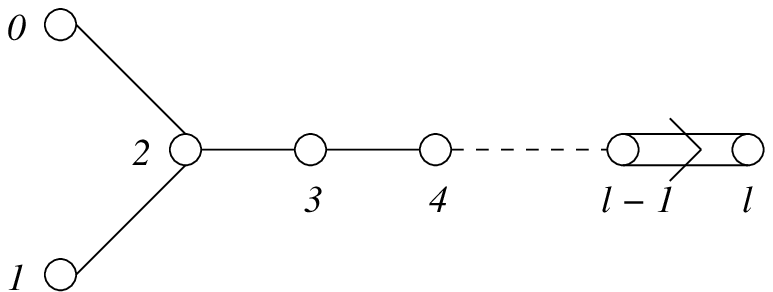}
}
\hfill}
}

There is an automorphism $\mu_2$ of the Dynkin diagram of $D_l^{(1)}$ obtained
by exchanging $0$ with $1$, exchanging $l-1$ with $l$, and fixing each of 
the other vertices.
The full heaps over $D_l^{(1)}$ discussed above fold via $\mu_2$ to a single
heap over the Dynkin
diagram of type $D_{l-1}^{(2)}$ (see Figure 14).  
In this case, heap elements formerly labelled $0$ or $1$ are relabelled by 
$0$, heap  elements formerly labelled $l-1$ or $l$ are relabelled by $l-2$,
and other heap elements have their labels decreased by $1$.  

\topcaption{Figure 14} The Dynkin diagram of type $D_{l+1}^{(2)}$
\endcaption
\centerline{
\hbox to 3.388in{
\vbox to 0.388in{\vfill
        \includegraphics{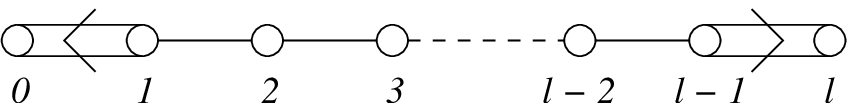}
}
\hfill}
}

\vfill\eject

\subhead Type $E_6^{(1)}$ \endsubhead

The Dynkin diagram $\Gamma$ of type $E_6^{(1)}$ is labelled as in Figure 15.

\topcaption{Figure 15} The Dynkin diagram of type $E_6^{(1)}$
\endcaption
\centerline{
\hbox to 2.138in{
\vbox to 1.319in{\vfill
        \includegraphics{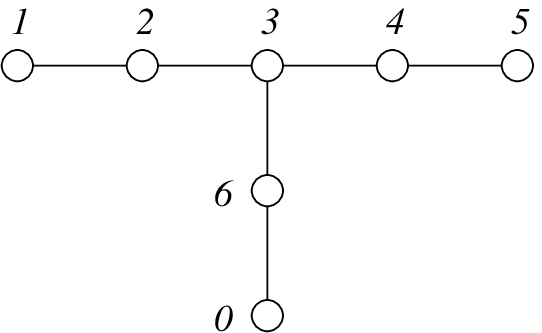}
}
\hfill}
}

There are two full heaps over $\Gamma$: the one shown
in Figure 16, and its dual, which may be constructed by applying by twisting
by a diagram automorphism corresponding to an odd
permutation to the branches of the Dynkin diagram emerging from vertex $3$.
The corresponding finite heaps are $F(E_6, 1)$ and $F(E_6, 5)$ in Wildberger's
notation.  There are $27$ orbits of proper ideals of $E$ under $\phi$.

\topcaption{Figure 16} A full heap over the Dynkin diagram of type $E_6^{(1)}$
\endcaption
\centerline{
\hbox to 1.527in{
\vbox to 2.527in{\vfill
        \includegraphics{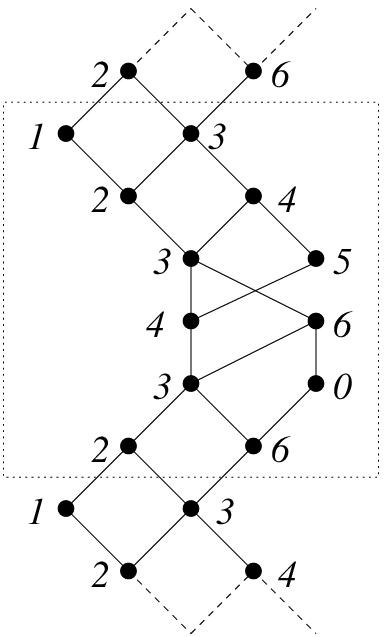}
}
\hfill}
}

\vfill\eject

\subhead Type $E_7^{(1)}$ \endsubhead

The Dynkin diagram $\Gamma$ of type $E_7^{(1)}$ is labelled as in Figure 17.

\topcaption{Figure 17} The Dynkin diagram of type $E_7^{(1)}$
\endcaption
\centerline{
\hbox to 3.138in{
\vbox to 0.819in{\vfill
        \includegraphics{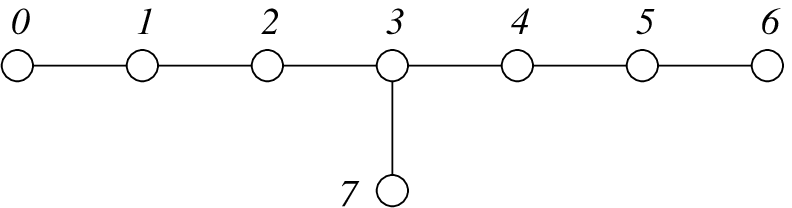}
}
\hfill}
}

There is one (self-dual) full heap over $\Gamma$, shown in Figure 18.  The
corresponding finite heap in Wildberger's notation is $F(E_7, 6)$.  There
are $56$ orbits of proper ideals of $E$ under $\phi$.

\topcaption{Figure 18} A full heap over the Dynkin diagram of type $E_7^{(1)}$
\endcaption
\centerline{
\hbox to 2.027in{
\vbox to 2.777in{\vfill
        \includegraphics{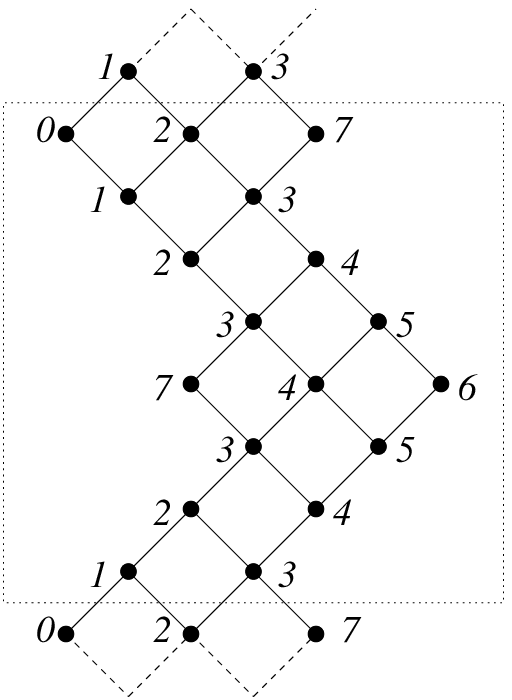}
}
\hfill}
}

\vfill\eject

\leftheadtext{} \rightheadtext{}

\end